# A Dynamic Bus Lane Strategy for Integrated Management of Human-Driven and Autonomous Vehicles

Haoran Li[a], Zhenzhou Yuan[a,*], Rui Yue[a,*], Guangchuan Yang[b], Fan Zhang[c], Zong Tian[d], Chuang Zhu[a]

[a] *School of Traffic and Transportation, Beijing Jiaotong University, 3 Shangyuancun, Haidian District, Beijing 100044, China;*

[b] *Institute of Transportation Research and Education, North Carolina State University, Raleigh, NC 27606, United States*

[c] *Institute of Remote Sensing and Geographic Information System, School of Earth and Space Sciences, Peking University, Beijing 100871, China*

[d] *University of Nevada, Reno 1664 N, Virginia Street, MS258, Reno, Nevada, 89557, United States*

*Corresponding authors: Rui Yue, Email: yuerui@bjtu.edu.cn

Zhenzhou Yuan, Email: zzy@bjtu.edu.cn

**Abstract:** This study introduces a dynamic bus lane (DBL) strategy, referred to as the dynamic bus priority lane (DBPL) strategy, designed for mixed traffic environments featuring both manual and automated vehicles. Unlike previous DBL strategies, this approach accounts for partially connected and autonomous vehicles (CAVs) capable of autonomous trajectory planning. By leveraging this capability, the strategy grants certain CAVs Right-of-Way (ROW) in bus lanes while utilizing their "leading effects" in general lanes to guide vehicle platoons through intersections, thereby indirectly influencing the trajectories of other vehicles. The ROW allocation is optimized using a mixed-integer linear programming (MILP) model, aimed at minimizing total vehicle travel time. Since different CAVs entering the bus lane affect other vehicles' travel times, the model incorporates lane change effects when estimating the states of CAVs, human-driven vehicles (HDVs), and connected autonomous buses (CABs) as they approach the stop bar. A dynamic control framework with a rolling horizon procedure is established to ensure precise execution of the ROW optimization under varying traffic conditions. Simulation experiments across two scenarios assess the performance of the proposed DBPL strategy at different CAV market penetration rates (MPRs). Results show that with a 20% MPR, the average travel time for private cars reduces by 14%, and at a 40% MPR, it decreases by 20%, while travel time for buses increases by no more than 0.6 seconds. Additionally, a sensitivity analysis quantifies the impacts of key parameters such as private car demand, bus arrival intervals, bus stop locations, and right-turn ratios.

## 1. Introduction

Developing strategies to prioritize buses is crucial for enhancing their operational efficiency and boosting the attractiveness of public transportation (Li et al. 2021; Zeng et al. 2020). Among these strategies, the implementation of bus lanes remains a primary approach to achieving spatial separation between buses and other vehicles. Extensive research has examined the operational performance and influential factors of bus lanes, highlighting their effectiveness in improveing bus priority (Cui et al. 2019; Farid et al. 2015; Shen et al. 2019). However, while bus lanes enhance the speed and reliability of bus services, they occupy road space originally designated for private cars. This reallocation of space can potentially exacerbating congestion on general lanes (Bai et al. 2017; Russo et al. 2022).

To address this issue, the concept of dynamic bus lanes (DBL) has been introduced, allowing private cars to use bus lanes when buses are not present (Dadashzadeh and Ergun 2018). Early research efforts, such as the intermittent bus lane proposed by Viegas et al. (2004), prohibited cars from entering the bus lane ahead of a bus and coordinated with modified traffic signals to clear cars already traveling in the bus lane. Eichler et al. (2006) refined this strategy by proposing the Bus Lane with Intermittent Priority (BLIP), requiring cars ahead of a bus to exit the bus lane without relying on bus signal priority. With advancements in wireless communication technologies, vehicles and infrastructure can communicate in the emerging connected environment, providing new opportunities for bus lane control. Wu et al (2018) introduced the bus lanes with intermittent and dynamic priority (BLIDP) strategy, which requires only cars within a designated "Clear Distance" ahead of a bus to exit the bus lane, while cars beyond this distance or trailing behind the bus were permitted to remain. Luo et al. (2022) extended this concept by proposing a dynamic bus lane with a moving block based on the cellular automata model, which dynamically adjusting the "Clear Distance" in a fully connected environment. Othman et al. (2022) compared the BLIDP strategy with the exclusive bus lane (EBL) strategy and mixed lane strategy in a fully connected vehicle environment, analyzing the conditions under which the BLIDP strategy is applicable. Xie (2022) proposed a cooperative dynamic bus lane system that introduced a “dynamic reserved zone” in front of a moving bus. This system requires cars in the bus lane that could impede bus operations to promptly vacate the reserved zone, aligning closely with the principles of the BLIDP strategy.

Previous research has primarily focused on rule-driven control strategies, which are applicable only under specific conditions, such as moderate traffic volumes (Othman et al. 2022). The emergence of highly automated or autonomous vehicles has enabled precise vehicle control, leading to a paradigm shift in DBL research (Ding et al. 2022; Hu et al. 2021). Unlike human-driven vehicles (HDVs), which stop and queue at stop bars when encountering red lights, connected and automated vehicles (CAVs) can avoid such stops through advanced trajectory planning (Yao and Li 2021; Yu et al. 2019). Moreover, the development of autonomous buses has also become a key aspect in the deployment of autonomous vehicles (Zhang et al. 2023).

Therefore, some studies have explored the use of bus lanes from the perspective of individual vehicle control, referred to in this paper as trajectory planning-driven strategies. For buses, researchers have proposed trajectory planning methods to improve the efficiency of connected autonomous buses (CABs) operating in bus lanes (Hu et al. 2021; Zhang et al. 2022). For private cars, Zhang et al. (2021) introduced a trajectory-based control method for bus priority lane that allows CAVs to use the bus lane in a fully connected and automated traffic environment, which was tested at a two-lane signalized intersection. Similarly, Chen et al. (2022) leveraged a rhythmic control to propose a novel scheme aimed at providing exclusive Right-of-Way (ROW) for buses while minimizing the combined travel costs of buses and private cars in a fully connected and automated environment. Wu et al. (2024) developed an automated intersection bus priority control model (AIM-BP) for fully connected and automated traffic, which introduces dynamic bus lanes to clear vehicles ahead of buses. Shan et al. (2024) proposed a trajectory planning method that allows CAVs to use the bus lane for lane changes without interfering with bus trajectories, and demonstrated its effectiveness through a case study of a road with one bus lane and one general lane.

In summary, previous research has made significant progress in improving bus lane utilization, but certain limitations remain. First, existing strategies have not effectively integrated bus lane usage rules with the trajectory planning capabilities of CAVs; these two approaches are still treated as separate strategies. Rule-driven control strategies manage bus lanes on a lane-wide basis through predefined rules but fail to account for environments where CAVs capable of trajectory planning. Conversely, while trajectory planning-driven strategies enable CAVs to use bus lanes through precise trajectory planning, they focus on optimizing individual vehicle movements rather than maximizing overall lane efficiency. For example, granting CAVs unrestricted access to bus lanes may diminish their distinct "leading effect" in general lanes (Ma et al. 2022). Therefore, a critical question to be addressed is "Which CAVs should be eligible to use the bus lane?" Answering this question is essential for achieving optimal operational efficiency across both the general and bus lanes.

Second, existing strategies primarily focus on straight-moving vehicles using bus lanes, with limited attention to right-turning vehicles. However, at intersections where right-turning demands are unavoidable, it becomes essential to allocate right-turn pockets adjacent to the bus lane to accommodate these turns while ensuring smooth bus operations (Dadashzadeh and Ergun 2018).

Therefore, the objective of this study is to propose a dynamic strategy for bus lane usage that:

(1) Functions effectively in a mixed manual-automated traffic environment with CAVs capable of trajectory planning.

(2) Selectively allows CAVs to use bus lanes while leveraging their leading effect to achieve system-wide operational efficiency across both general and bus lanes.

(3) Accommodates the demands of right-turning vehicles and ensures the strategy's applicability in complex road scenarios.

To achieve these objectives, this study introduces a novel control strategy, the Dynamic Bus Priority Lane (DBPL), designed for mixed traffic environments involving HDVs, CAVs, and CABs, with the latter two capable of self-planning their trajectories. The DBPL strategy allocates ROW to selected CAVs by issuing lane-changing recommendations, optimizing the use of time-space resources in bus lanes. It establishes customized driving protocols to manage different vehicle types with distinct turning requirements when accessing bus lanes.

Furthermore, this research expanded the concept of using virtual vehicles for estimating the passing states of vehicle by accounting for roads with right-turn pockets, which improves the applicability of this method (Li et al. 2024). The optimization of ROW for CAVs is formulated as a Mixed Integer Linear Programming (MILP) problem (Zhu et al. 2024), leveraging the leading effect of CAVs to reduce start-up delays in mixed traffic. The model is efficiently executed using a rolling horizon control scheme that dynamically updates vehicle status in response to real-time traffic conditions. Additionally, a ROW Pre-allocation Heuristic Algorithm is proposed to improve solution efficiency by reducing the model's dimensionality.

The remainder of this paper is structured as follows: Section 2 outlines the research problems addressed in this study. Section 3 formulates the MLIP model for optimizing ROW allocation and introduces the solution method. Section 4 presents a series of simulation experiments to validate the performance of DBPL, comparing it with EBL control strategies under various parameters. Lastly, Section 5 summarizes concluding remarks and suggests directions for future research.

## 2. Problem statement

As illustrated in Fig. 1, we consider a signalized intersection under a mixed manual-automated traffic environment. A no-changing zone near the stop bar prohibits lane changes, requiring vehicles to complete lane changes before entering this area. The lanes are divided into two types: general lanes, accommodating both CAVs and HDVs, and a designated bus lane for CABs. CAVs and HDVs are treated as homogeneous in terms of physical dimensions (e.g., length and width), maximum speed, and maximum acceleration and deceleration capabilities. Within the intersection's communication area, the control center can obtain real-time information such as vehicle speed and position. The task of the control center is to optimize the ROW and the timing of CAVs entering the bus lane within a predefined planning horizon, based on both current and projected vehicle motion states. The objective is to minimize the weighted travel times of both cars and buses.

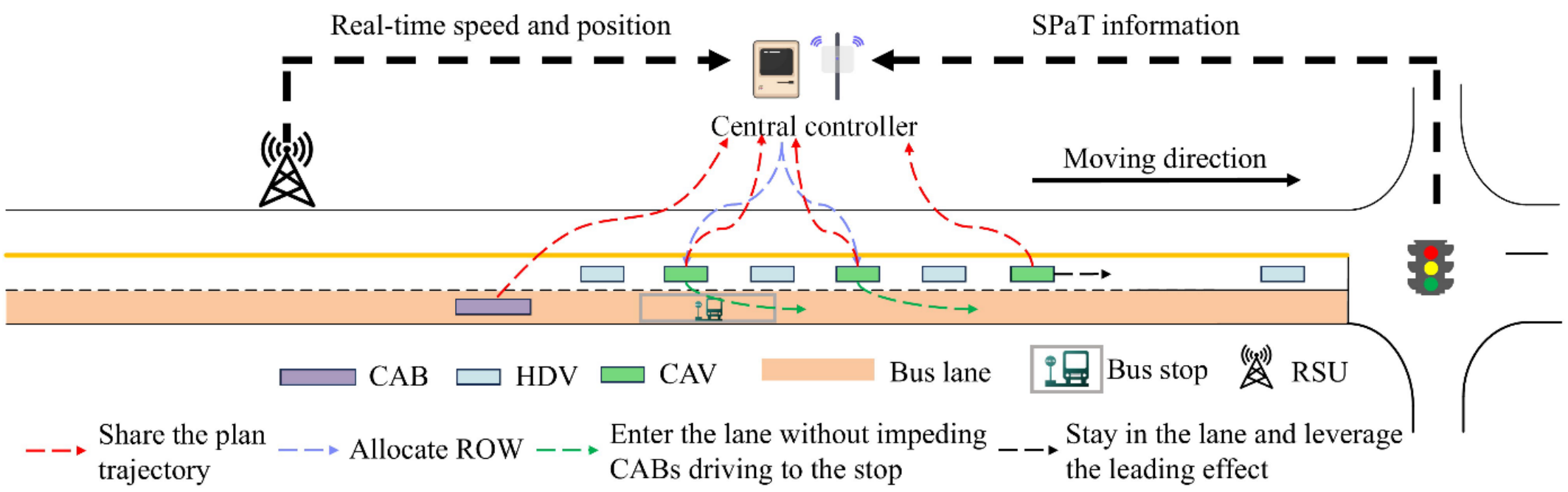

Fig. 1 The proposed DBPL control strategy.

CAVs and CABs act as intelligent agents that independently plan their trajectories using data on vehicle states and signal timings. To enter the bus lane, CAVs must share their planned trajectories with the control center for system coordination. CAVs granted ROW execute lane changes and replan their trajectories at designated times. CABs, when planning their trajectories, must consider three key processes: approaching the bus stop, dwelling at the stop, and passing the stop bar. They also need to account for the combined effects of bus stops and signalized intersections (Shan et al. 2023). Unlike the strategy proposed by Shan et al. (2024), which treats CAB trajectories as fixed constraints for planning CAV trajectories in the bus lane, this study ensures CAB priority through weighted coefficients in the objective function of the optimization. As a result, CAB priority may be slightly adjusted, and trade-offs will be analyzed in Section 4.

Since lane changes can impact the timing of subsequent vehicles passing the stop bar (as shown in Fig. 2, illustrating scenarios where lane changes may or may not affect subsequent vehicles), it is crucial to assess how these lane changes influence others to identify the optimal vehicles for entering or exiting the bus lane. The optimization model incorporates the distinct driving behaviors of HDVs and CAVs, with the leading effect of CAVs naturally emerging within the model.

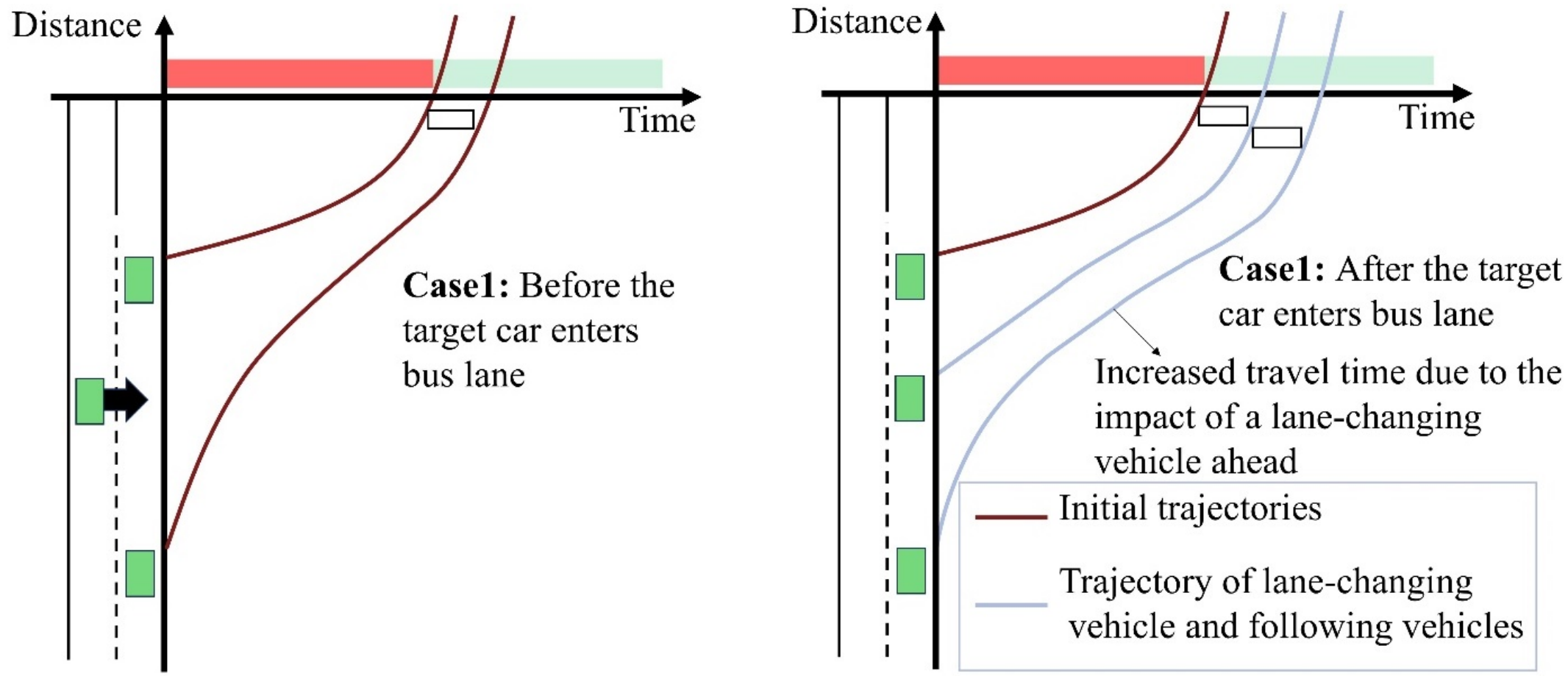

(a) Following vehicle is affected by the lane changing of the preceding vehicle

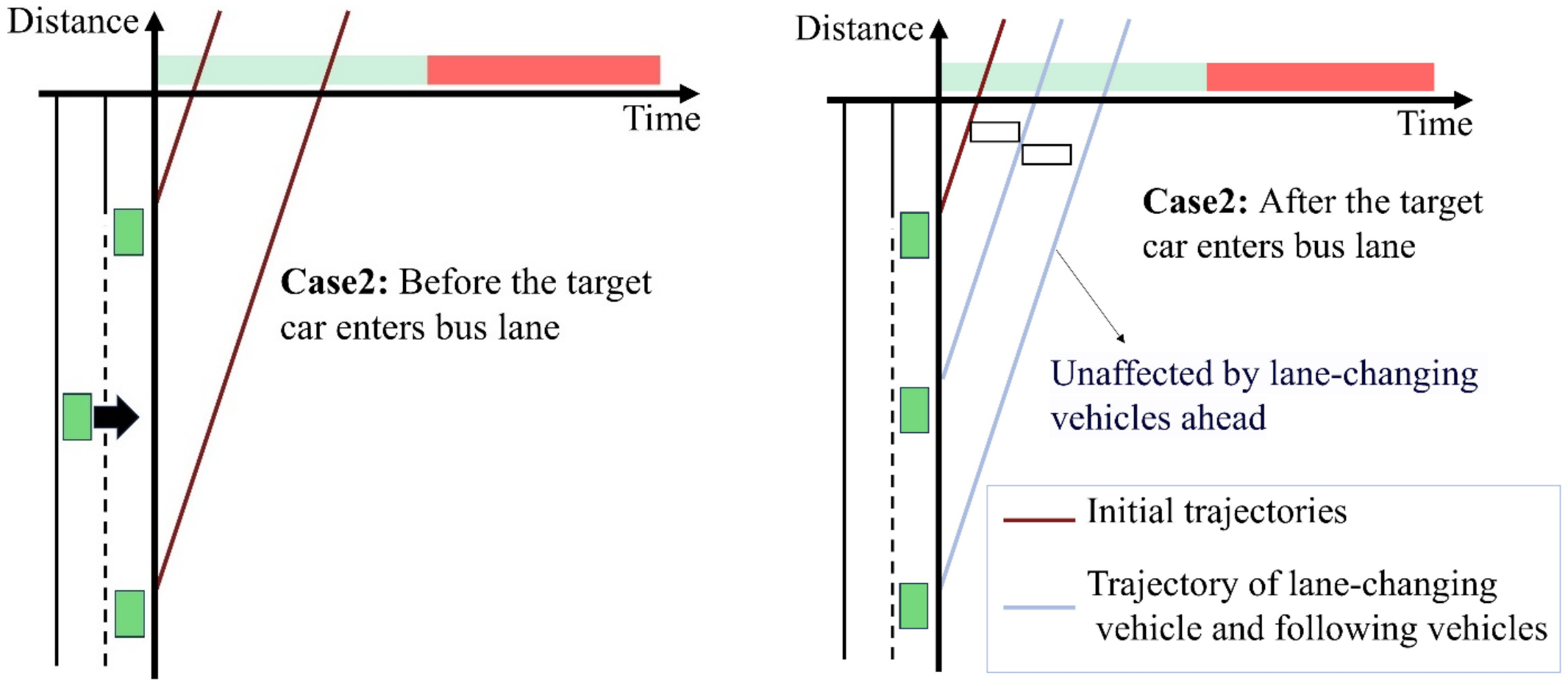

(b) Following vehicle is not affected by the lane changing of the preceding vehicle

Fig. 2 Schematic diagram of the influence of lane change on the following vehicle

# 3. Methodology

This section formulates a MILP model for optimizing ROW allocation. The objective of the model is to minimize the weighted average travel costs for both buses and cars, thereby enhancing the overall efficiency of the system. Specifically, the average travel time for cars, denoted as $T_c$, includes both travel time in the general lane and the bus lane. It is expressed as follows:

$$T_c = \frac{1}{|\mathcal{I}_h \cup \mathcal{I}_a|} \cdot \sum_{k \in K} \sum_{i \in \mathcal{I}_h \cup \mathcal{I}_a} \lambda(k) \cdot t_i^{dep}(k) \tag{1}$$

where $i$ is the vehicle index. $t_i^{dep}(k)$ represents the estimated time for vehicle $i$ to pass the stop bar at time $k$. $K$ is the set of all time steps within the planning horizon: $K = \{k_0, k_0 + \Delta k, k_0 + 2\Delta k, \ldots, k_0 + h\}$, where $k_0$ is the start time step and $h$ is the length of the planning horizon. $\lambda(k)$ is binary variable, with $\lambda(k) = 1$ indicating that certain CAVs are recommended to change lanes at time $k$, and $\lambda(k) = 0$ otherwise. $\mathcal{I}_h$ and $\mathcal{I}_a$ are the sets of HDVs and CAVs within the study area, respectively. $|\mathcal{I}_h \cup \mathcal{I}_a|$ denotes the total number of vehicles in the combined set of HDVs and CAVs.

Similarly, the average travel time for buses, $T_b$, is calculated by:

$$T_b = \frac{1}{|\mathcal{I}_b|} \cdot \sum_{k \in K} \sum_{i \in \mathcal{I}_b} \lambda(k) \cdot t_i^{dep}(k) \tag{2}$$

where $\mathcal{I}_b$ represents the set of CABs within the study area. To prioritize buses, a weighting parameter $\omega_p$, ranging from [0,1], is introduced to balance the efficiency priority between buses and private cars. The objective function is then formulated as:

$$\min Z = \omega_p \cdot T_b + (1 - \omega_p) \cdot T_c \tag{3}$$

The average travel time of vehicles in the objective function Eq.(3) can be calculated based on vehicle motion estimation from Section 3.1. The future motion of vehicle $i$ depends on its current state, the leading vehicle, and the signal timing information. The vehicle states (i.e., location $x_i(k)$ and speed $v_i(k)$) are denoted as:

$$S = \{x_i(k), v_i(k) | i \in \mathcal{I}_b \cup \mathcal{I}_a \cup \mathcal{I}_h, k \in K\} \tag{4}$$

The future states of CAVs and CABs within the range $[k_0 + \Delta k, k_0 + h]$ can be derived from shared planned trajectories, while the future states of HDVs in the bus lane during the same time range are predicted using the Intelligent Driver Model (IDM).

The decision variable determines whether each vehicle in the general lane can change into the bus lane at time $k$. The set of these decision variables is denoted as $\Phi$:

$$\Phi = \{\phi_i(k) | i \in \mathcal{I}_{gl}, k \in K\} \tag{5}$$

where $\mathcal{I}_{gl}$ represents the set of vehicles in the general lane at current time $k_0$, and $\phi_i(k)$ is a binary variable. If vehicle $i$ changes lanes at time $k$, then $\phi_i(k) = 1$; otherwise, $\phi_i(k) = 0$.

Signal timings, which can be obtained in advance, are denoted as $Y$. Therefore, the travel cost $Z$ can be expressed as a function of $S$, $Y$, and $\Phi$:

$$Z = f(S, Y, \Phi) \tag{6}$$

The relationship between $Z$ and $\Phi$, given as $S$ and $Y$, is further detailed in Sections 3.1 and 3.2.

### *3.1 Vehicle lateral motion modeling*

This section focuses on estimating vehicle travel time when individual lane-changing behaviors are uncertain, given the vehicle states $S$ and signal timings $Y$. Unlike previous studies that primarily examined longitudinal driving behavior, this study also considers lateral vehicle movements, as some CAVs in the general lane may switch to the bus lane. Typically, a vehicle's future motion is influenced by the signal timings of its current lane and the motion of preceding vehicles. If a preceding vehicle changes lanes, the motion of the following vehicle may also be affected. However, when the final vehicle distribution in the lane is unknown, estimating their passing states becomes challenging.

Therefore, this study introduces a modeling approach using "virtual vehicles" (Li et al. 2024) to estimate each vehicle's new travel time based on $S$ and $Y$, without prior knowledge of which CAVs will change lanes. Specifically, any CAV $i$ in the general lane is mapped to a virtual vehicle $r_i$ in the bus lane. The set of virtual vehicles in the bus lane is denoted as $\mathcal{I}_{bl}^{vir}$, and the set of real vehicles as $\mathcal{I}_{bl}^{real}$, with $\mathcal{I}_{bl} = \mathcal{I}_{bl}^{vir} \cup \mathcal{I}_{bl}^{real}$. The characteristics of virtual vehicles are defined as follows:

(1) When CAV $i$ does not change lanes, the $v_{r_i}(k)$, $x_{r_i}(k)$ and $t_{r_i}^{dep}(k)$ of the virtual vehicle $r_i$ inherit those of its preceding vehicle, such as vehicle $r_{i_1}$ in Fig. 3.

(2) When CAV $i$ changes lanes, the virtual vehicle $r_i$ is converted into a real vehicle, inheriting the speed $v_{r_i}(k)$ and position $x_{r_i}(k)$ of CAV $i$, and its new departure time $t_{r_i}^{dep}(k)$ is recalculated, as illustrated by vehicle $r_{i_2}$ in Fig. 3. Simultaneously, the real vehicle in the general lane is converted into a virtual vehicle, inheriting the properties of its preceding vehicle, such as vehicle $i_2$ in Fig. 3. This method allows the estimation of passing states when vehicles change lanes and accurately reflects the new vehicle-following dynamics in both lanes.

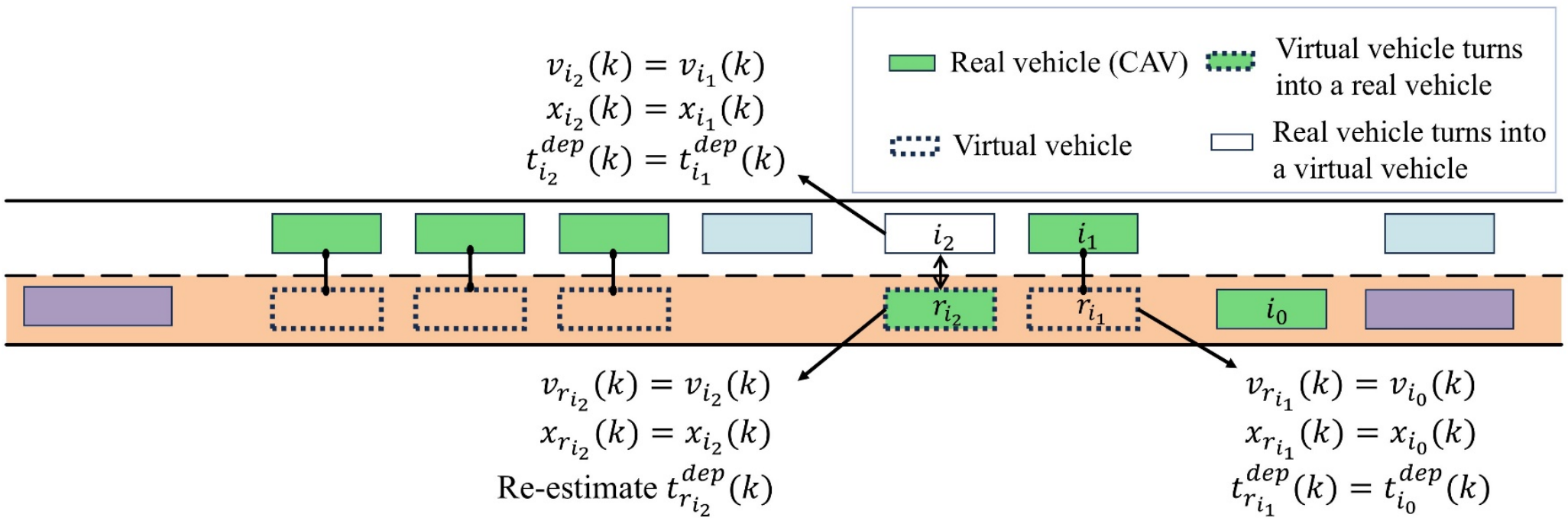

Fig. 3 Schematic diagram of virtual vehicles modeling method

Using the virtual vehicles concept, $T_c$ is reformulated as follows:

$$T_c = \frac{\lambda(k)}{|\mathcal{I}_h \cup \mathcal{I}_a|} \cdot \left( \sum_{k\in K} \sum_{i\in \mathcal{I}_{gl}} (1-\phi_i(k)) \cdot t_i^{dep}(k) + \sum_{k\in K} \sum_{i\in \mathcal{I}_{bl}^{real} \setminus \mathcal{I}_b} t_i^{dep}(k) + \sum_{k\in K} \sum_{i\in \mathcal{I}_{bl}^{vir}} \phi_i(k) \cdot t_i^{dep}(k) \right) \tag{7}$$

The first term represents the passing time of private cars in set $\mathcal{I}_{gl}$ that remain their current lane. The second term refers to the passing time of private cars already in the bus lane. The third term captures the passing time of private cars selected to enter the bus lane (i.e., changing from virtual to real vehicles).

For CAV $i$ in the general lane, it can change into the bus lane at most once during the planning horizon:

$$\sum_{k\in K} \phi_i(k) \le 1, \forall i \in \mathcal{I}_{gl} \cap \mathcal{I}_a \tag{8}$$

Moreover, the control center sends a lane-change recommendation to a CAV only once within the planning horizon.

$$\sum_{k\in K} \lambda(k) = 1 \tag{9}$$

Only CAVs that satisfy lateral operation constraints are eligible for lane changes. First, CAVs cannot change lanes when approaching the no-changing zone near intersections:

$$x_i(k) - M \cdot (1-\phi_i(k)) \le x_c - x_n, \forall i \in \mathcal{I}_{gl} \cap \mathcal{I}_a; k \in K \tag{10}$$

where $x_c$ is the stop bar position; and $x_n$ is the starting position of the no-changing zone.

Second, lane changes are not allowed if CAV $i$ is stopped:

$$-M \cdot v_i(k) \le \phi_i(k) \le M \cdot v_i(k), \forall i \in \mathcal{I}_{gl} \cap \mathcal{I}_a; k \in K \tag{11}$$

For safety, CAVs must meet lateral safety constraints. This model assumes CAVs can complete lane changes within one time step. If the control center recommends a lane change, there must be enough space in the bus lane to avoid collisions:

$$\left(1-\phi_i(k)\right)\cdot M + x_{\overline{r}_i}(k) - x_i(k) > d_{safe}, \forall i \in \mathcal{I}_{gl} \cap \mathcal{I}_a; k \in K \tag{12}$$

$$\left(1-\phi_i(k)\right)\cdot M + x_i(k) - x_{\underline{r}_i}(k) > d_{safe}, \forall i \in \mathcal{I}_{gl} \cap \mathcal{I}_a; k \in K \tag{13}$$

where $d_{safe}$ is the safe distance that the vehicle needs to maintain with the preceding vehicle and the following vehicle in the target lane for lane-changing. $\overline{r}_i$ and $\underline{r}_i$ are the preceding vehicle and the following vehicle of vehicle $i$ in the target lane, respectively.

Upon completing a lane change, the vehicle should maintain a safe distance from the preceding vehicle $\overline{r}_i$ and the following vehicle $\underline{r}_i$ to ensure no collision occurs after the lane change, where $k_{lc}$ represents the lane changing time step:

$$x_i(k+k_{lc}) + M\cdot(1-\phi_i(k)) < x_{\overline{r}_i}(k+k_{lc}), \forall i \in \mathcal{I}_{gl} \cap \mathcal{I}_a; k \in K \tag{14}$$

$$x_i(k+k_{lc}) + M\cdot(1-\phi_i(k)) > x_{\underline{r}_i}(k+k_{lc}), \forall i \in \mathcal{I}_{gl} \cap \mathcal{I}_a; k \in K \tag{15}$$

### *3.2 Vehicle longitudinal motion modeling*

In this paper, Newell's car-following model is used to capture the influence of a preceding vehicle $\bar{i}$ on the passing state of the subject vehicle $i$ being estimated (Newell 2002). The longitudinal vehicle passing state estimation assumes that vehicles aim to traverse the stop bar smoothly as soon as possible (Ma et al. 2023; Sharma et al. 2018). Building on the longitudinal motion estimation method and the concept of virtual vehicles, this paper further incorporates the impact of lane changes on vehicle passing state estimation. Additionally, the analysis includes right-turning vehicles transitioning into the designated right-turn pocket, as illustrated in Fig. 4. Through lane-based motion analysis, we determine the estimated time $t_i^{dep}(k)$ for each vehicle to pass the stop bar at each time step $k$ within the planning horizon.

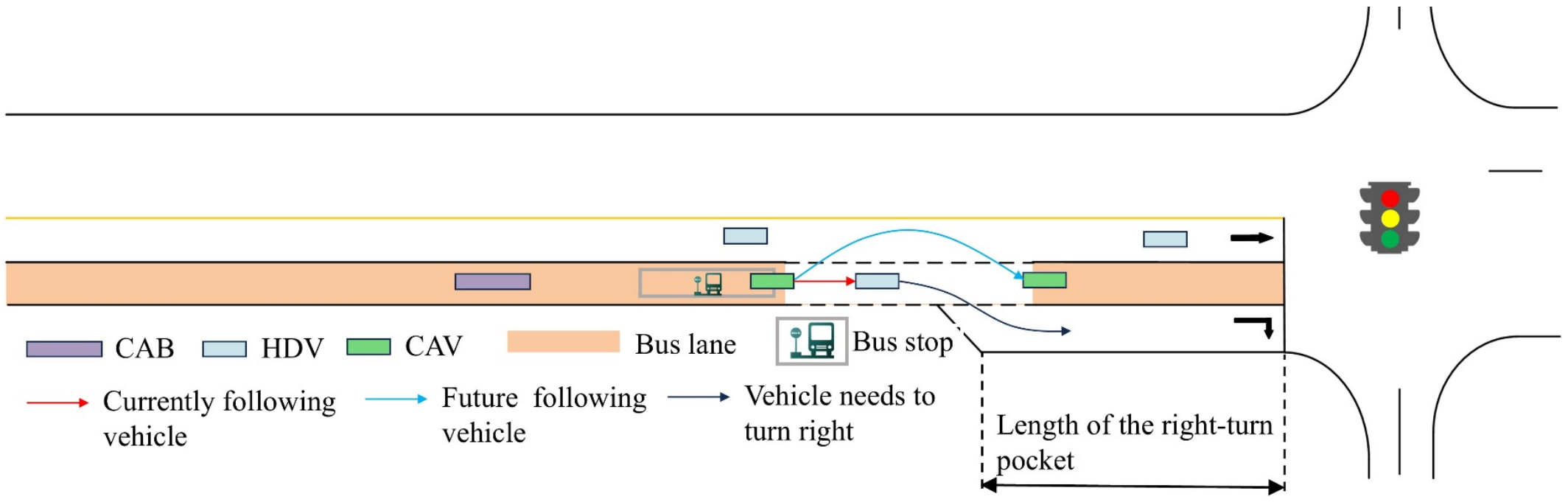

Fig. 4 Roadway with right-turn pocket

Within the planning horizon $K$, vehicle speed and position vary with each time step $k$. To accurately estimate passing states, the original vehicle sets are divided into different subsets.

(1) For vehicles in the bus lane ($\mathcal{I}_{bl}$)

- $I'_{bl}(k)$: Vehicles that have passed the stop bar at time step $k$:

$$I'_{bl}(k) = \{i | x_i(k) \geq x_c, i \in \mathcal{I}_{bl}\} \tag{16}$$

where $x_c$ is the position of the stop bar. Since right-turning vehicles in the bus lane need to change lanes upon reaching the right-turn pocket entrance, the set $\mathcal{I}_{bl}$ is further divided at time step $k$ into:

- $\overline{I}_{bl}(k)$: Vehicles between the right-turn pocket entrance and the stop bar:

$$\overline{I}_{bl}(k) = \{i | x_w \leq x_i(k) \leq x_c, i \in \mathcal{I}_{bl}\} \tag{17}$$

- $\underline{I}_{bl}(k)$: Vehicles from the entrance of the roadway to the right-turn pocket entrance:

$$\underline{I}_{bl}(k) = \{i | 0 \leq x_i(k) \leq x_w, i \in \mathcal{I}_{bl}\} \tag{18}$$

where $x_w$ is the postion of right-turn pocket entrance. Fig. 5 provides an example of how $\mathcal{I}_{bl}$ is divided into subsets $I'_{bl}(k)$, $\overline{I}_{bl}(k)$, and $\underline{I}_{bl}(k)$. At the initial time step $k_0$, the set $\mathcal{I}_{bl} = \{i_0, i_1, i_2, i_3, i_4, i_5\}$. At this point, $I'_{bl}(k_0) = \emptyset$, $\overline{I}_{bl}(k_0) = \emptyset$, $\underline{I}_{bl}(k_0) = \{i_0, i_1, i_2, i_3, i_4, i_5\}$. As time progresses, the composition of these subsets changes, as shown in time steps $k_3$ and $k_6$ in Fig. 5. Since $\underline{I}_{bl}(k)$ includes both right-turning and through vehicles, the lane-changing behavior of the preceding right-turning vehicles needs to be considered when estimating the passing states of through vehicle. For instance, in Fig. 5, vehicle $i_5$ initially follows the right-turning vehicle $i_4$. Once $i_4$ enters the right-turn pocket, the preceding vehicle for $i_5$ becomes $i_3$. This change in the following relationship is addressed through the modeling described in Section 3.2.1.

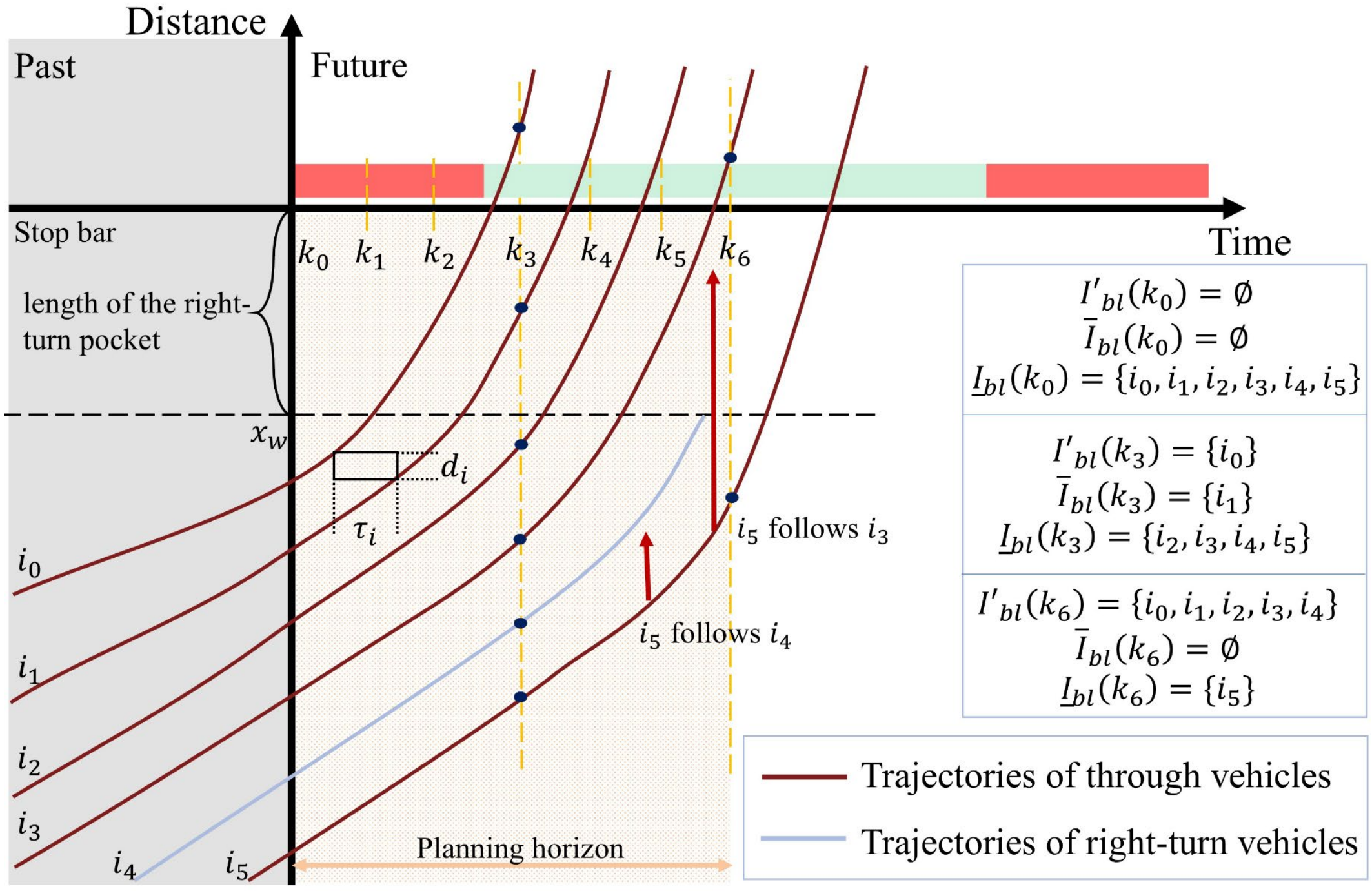

Fig. 5 Schematic diagram of vehicle subset division method in bus lane

(2) For vehicles in the general lane ($\mathcal{I}_{gl}$)

- $I'_{gl}(k)$: Vehicles that have passed the stop bar at time step $k$:

$$I'_{gl}(k) = \{i | x_i(k) \geq x_c, i \in \mathcal{I}_{gl}\} \tag{19}$$

- $\bar{I}_{gl}(k)$: Vehicles between the entrance of the roadway and the stop bar at time step $k$:

$$\bar{I}_{gl}(k) = \{i | 0 \leq x_i(k) \leq x_c, i \in \mathcal{I}_{gl}\} \tag{20}$$

The vehicle passing state estimation process is illustrated in Fig. 6.

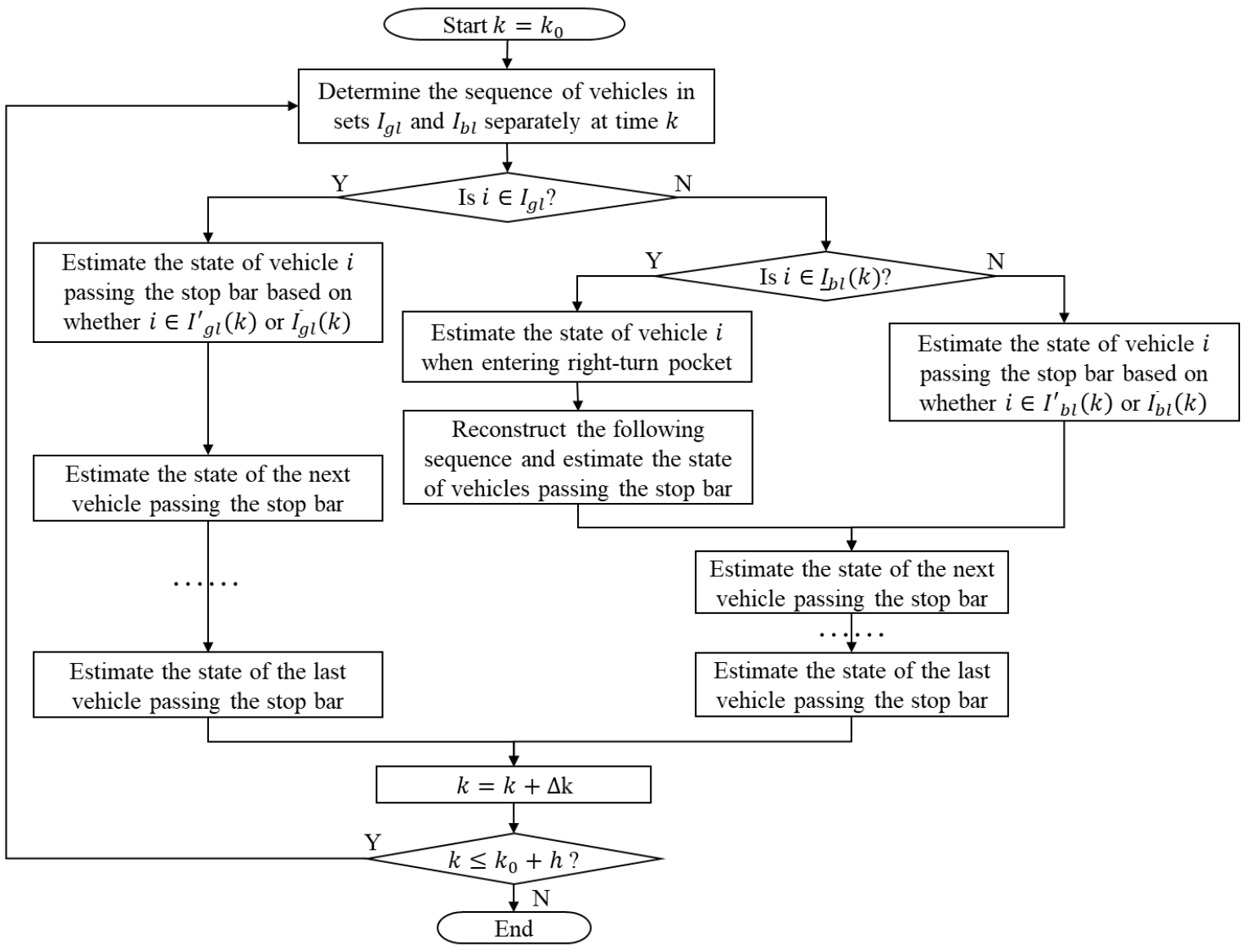

Fig. 6 The scheme of vehicle passing states estimation

### *3.2.1 Estimation of passing states for vehicle in bus lane*

(1) Estimation of the departure time $t_i^{dep}(k)$ for vehicles in set $I'_{bl}(k)$:

Within the planning horizon $K$, if vehicle $i$ in $I'_{bl}(k)$ is a real vehicle, $t_i^{dep}(k)$ can be obtained through real-time shared trajectory data, denoted as $t_i^{tra}(k)$. If it is a virtual vehicle, its $t_i^{dep}(k)$ is $t_{bl}$ (the time the last real vehicle passed the stop bar during the previous green phase) or the departure time $t_{\bar{i}}^{dep}(k)$ of the preceding vehicle. The formula for $t_i^{dep}(k)$ is:

$$t_i^{dep}(k) = \begin{Bmatrix} t_i^{tra}(k), if\ i \in I'_{bl}(k) \cap \mathcal{I}_{bl}^{real} \\ t_{bl}, if\ i\ is\ the\ first\ vehicle\ in\ I'_{bl}\ and\ i \in I'_{bl}(k) \cap \mathcal{I}_{bl}^{vir} \\ t_{\bar{i}}^{dep}(k), if\ i\ is\ not\ the\ first\ vehicle\ in\ I'_{bl}\ and\ i \in I'_{bl}(k) \cap \mathcal{I}_{bl}^{vir} \end{Bmatrix}, \forall k \in K \tag{21}$$

(2) Estimation of the departure time $t_i^{dep}(k)$ for vehicles in set $\bar{I}_{bl}(k)$:

Assuming that vehicles immediately enter the right-turn pocket upon encountering it, all vehicles in the set $\bar{I}_{bl}(k)$ are treated as through vehicles, as right-turning vehicles will not appear in this set. We denote $t_i^{p1}(k)$ as the earliest time at which a vehicle can pass the stop bar, without considering its preceding vehicle or signal control. The formula for $t_i^{p1}(k)$ is:

$$t_i^{p1}(k) = k + Z(x_c - x_i(k); v_i(k); v_i^{max}; a_i^{max}), \forall i \in \bar{I}_{bl}(k) \cup \bar{I}_{gl}(k), k \in K \tag{22}$$

where $v_i^{max}$ and $a_i^{max}$ denote the maximum speed and acceleration of vehicle $i$, respectively. The function $Z(\cdot)$ calculates the minimum time required to travel a given distance, considering the vehicle's initial speed, maximum speed, and acceleration limits (Zheng et al. 2020). The relevant formulas are:

$$t_i^{acc}(k) = \frac{v_i^{max} - v_i(k)}{a_i^{max}}, \forall i \in \bar{I}_{bl}(k) \cup \bar{I}_{gl}(k), k \in K \tag{23}$$

$$s_i^{acc}(k) = \frac{v_i^{max} + v_i(k)}{2} \cdot t_i^{acc}(k), \forall i \in \bar{I}_{bl}(k) \cup \bar{I}_{gl}(k), k \in K \tag{24}$$

$$t_i^{p1}(k) = \begin{Bmatrix} \frac{-v_i(k) + \sqrt{(v_i(k))^2 + 2a_i^{max}(x_c - x_i(k))}}{a_i^{max}}, \ if\ s_i^{acc}(k) > x_c - x_i(k) \\ \frac{x_c - x_i(k) - s_i^{acc}(k)}{v_i^{max}} + t_i^{acc}(k), otherwise \end{Bmatrix}, \forall i \in \bar{I}_{bl}(k) \cup \bar{I}_{gl}(k), k \in K \tag{25}$$

The time a vehicle passes the stop bar, accounting for the influence of preceding vehicles, is denoted as $t_i^{p2}(k)$:

$$t_i^{p2}(k) = max(t_i^{p1}(k), t_i^{pre}(k) + \hat{\tau}_i), \forall i \in \bar{I}_{bl}(k) \cup \bar{I}_{gl}(k), k \in K \tag{26}$$

where $\hat{\tau}_i$ is the minimum time headway, a simplified substitute for the function $Z(d_i; v_i(k); v_i^{max}; a_i^{max})$ (Ma et al. 2023). $t_i^{pre}(k)$ is the time when the preceding vehicle of $i$ passes the stop bar. Its value depends on the position of vehicle $i$ along the road:

- If vehicle $i$ is the first vehicle at both time $k_0$ and $k$: $t_i^{pre}(k) = t_{bl}$, as exemplified by vehicle $i_1$ at time $k_{n-1}$ in Fig. 7. Here, $t_{bl}$ is the time when the last actual vehicle $i_0$ passes the stop bar during the previous green light phase.
- If vehicle $i$ is initially a following vehicle at time $k_0$ but becomes the leading vehicle at time $k$: $t_i^{pre}(k) = t_{\bar{i}}^{tra}(k)$, as demonstrated by vehicle $i_3$ at time $k_n$ in Fig. 7. Here, $t_{i_2}^{tra}(k)$ is the time when vehicle $i_2$ passes the stop bar obtained through known vehicle states $S$.
- If vehicle $i$ is not the first vehicle at time $k$: $t_i^{pre}(k) = t_{\bar{i}}^{dep}(k)$, as exemplified by vehicle $i_3$ at time $k_{n-1}$ in Fig. 7.

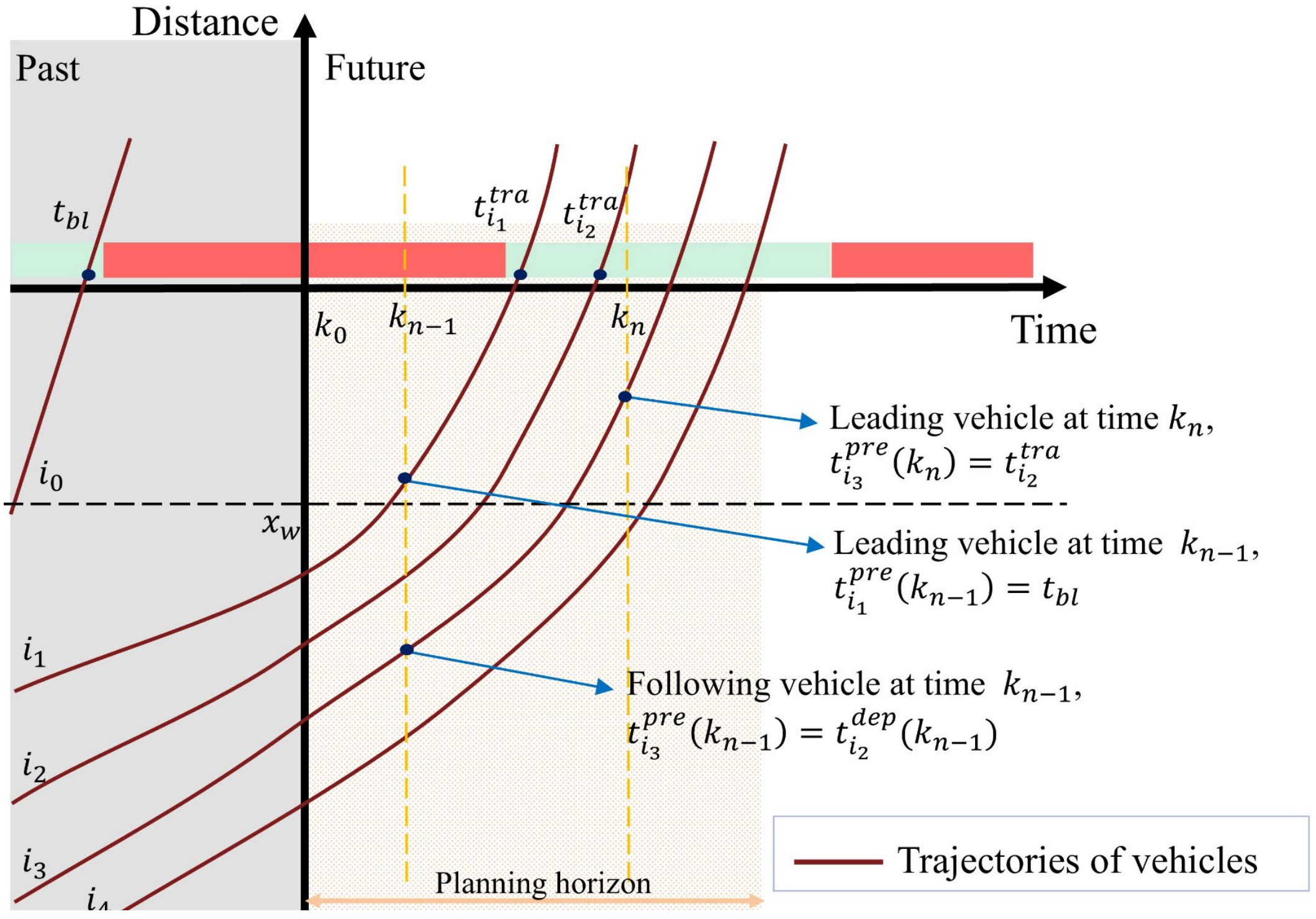

Fig. 7 Schematic diagram of the method for judging the value of $t_i^{pre}(k)$

When vehicle $i$ passes the stop bar, the signal timing also influences the departure time, calculated as:

$$t_i^{p3}(k) = max(t_i^{p2}(k), t_R + t_r), \forall i \in \bar{I}_{bl}(k) \cup \bar{I}_{gl}(k), k \in K \tag{27}$$

where $t_r$ is the red phase duration, and $t_R$ is the starting time of the red phase within the signal cycle, determined as:

$$t_R = \left\lfloor \frac{t_i^{p2}(k)}{t_c} \right\rfloor \cdot t_c, \forall i \in \bar{I}_{bl}(k) \cup \bar{I}_{gl}(k), k \in K \tag{28}$$

where $t_c$ is the signal cycle length ( $t_c = t_r + t_g$).

Finally, the departure time $t_i^{dep}(k)$ is determined as:

$$t_i^{dep}(k) = \begin{cases} t_i^{p3}(k) + (1 - \lambda(k)) \cdot M, \ if \ i \in \bar{I}_{bl}(k) \cap \mathcal{I}_{bl}^{real} \\ \phi_{r_i}(k) \cdot t_i^{p3}(k) + (1 - \phi_{r_i}(k)) \cdot t_i^{pre}(k) + (1 - \lambda(k)) \cdot M, otherwise \end{cases}, \forall i \in \bar{I}_{bl}(k) \cup \bar{I}_{gl}(k), k \in K \tag{29}$$

(3) Estimation of the departure time $t_i^{dep}(k)$ vehicles in set $\underline{I}_{bl}(k)$:

In practice, bus lanes are typically positioned at the outermost edge of the road, leading to interactions between through vehicles and right-turning vehicles. Right-turning vehicles temporarily borrow the bus lane before entering the right-turn pocket. This scenario complicates the estimation of passing states within the bus lane compared to general lanes since the behavior of right-turning vehicles directly affects the passing states of following through vehicles (as illustrated in Fig. 8). For instance, through vehicles may either follow a right-turning vehicle for a certain distance, pass the stop bar unhindered (like $i_5$), or follow other vehicles (like $i_2$).

To accurately estimate the time for vehicles to pass the stop bar, we adopt a two-step approach:

Firstly, assess the passing states, including the passing time $\underline{t_i}^{dep}(k)$ and speeds $\underline{v_i}^{dep}(k)$ of each vehicle in set $\underline{I}_{bl}(t)$ as they pass $x_w$. Secondly, reorganizing the car-following relationships based on whether the vehicle intends to turn right or continue through and estimating the passing time $t_i^{dep}(k)$ for through vehicles.

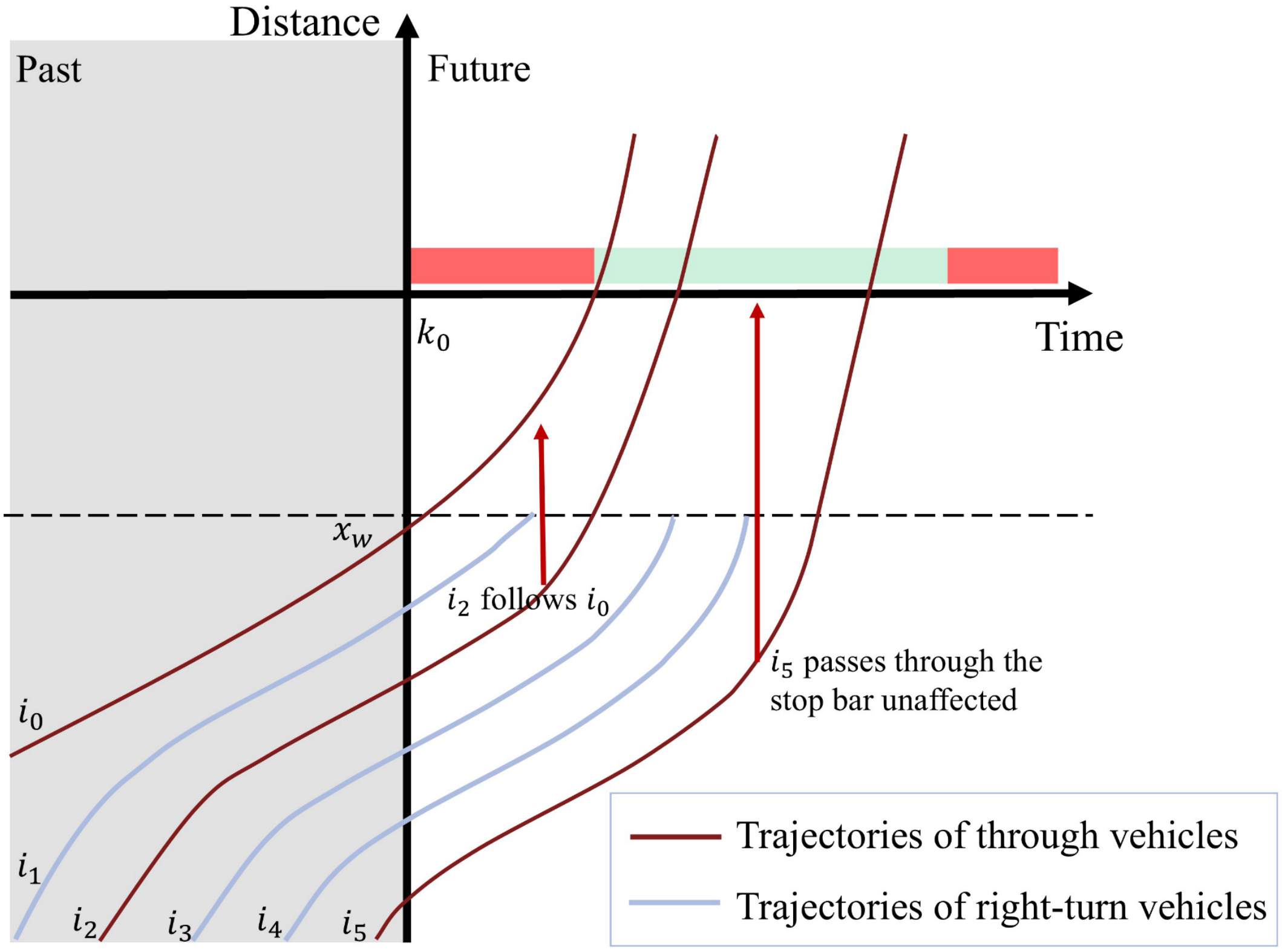

Fig. 8 Diagram of car-following relationship alteration

The earliest time a vehicle $i$ can pass $x_w$ without considering preceding vehicles is denoted as $\underline{t}_i^{p1}(k)$, calculated as:

$$\underline{t}_i^{p1}(k) = k + Z(x_w - x_i(k); v_i(k); v_i^{max}; a_i^{max}), \forall i \in \underline{I}_{bl}(k), k \in K \tag{30}$$

Taking into account the influence of preceding vehicles, the adjusted passing time $\underline{t}_i^{p2}(k)$ is given by:

$$\underline{t}_i^{p2}(k) = max(\underline{t}_i^{p1}(k), \underline{t}_i^{pre}(k) + \hat{\tau}_i), \forall i \in \underline{I}_{bl}(k), k \in K \tag{31}$$

Here, $\underline{t}_i^{pre}(k)$ is the time when the preceding vehicle $\bar{i}$ passes $x_w$, and $\underline{v}_i^{pre}(k)$ is the speed of $\bar{i}$ at that point. The values of $\underline{t}_i^{pre}(k)$ and $\underline{v}_i^{pre}(k)$ are determined as follows:

- If vehicle $i$ is not the first vehicle within the segment $[0, x_w]$ at time $k$, then $\underline{t}_i^{pre}(k) = \underline{t}_{\bar{i}}^{dep}(k)$ and $\underline{v}_i^{pre}(k) = \underline{v}_{\bar{i}}^{dep}(k)$.
- If vehicle $i$ is the first vehicle within $[0, x_w]$ both at time $k_0$ and $k$, then $\underline{t}_i^{pre}(k) = \underline{t}_{bl}$ and $\underline{v}_i^{pre}(k) = \underline{v}_{bl}$. Here, $\underline{t}_{bl}$ and $\underline{v}_{bl}$ are the time and speed when the last real vehicle passed $x_w$.
- If vehicle $i$ becomes the first vehicle within the segment $[0, x_w]$ at $k$, though it was not at $k_0$,then $\underline{t}_i^{pre}(k) = \underline{t}_{\bar{i}}^{tra}(k)$ and $\underline{v}_i^{pre}(k) = \underline{v}_{\bar{i}}^{tra}(k)$. $\underline{t}_{\bar{i}}^{tra}(k)$ and $\underline{v}_{\bar{i}}^{tra}(k)$ refer to the time and speed at which the preceding vehicle passed $x_w$, derived from the vehicle state set $S$.

The time for vehicle $i$ to pass $x_w$ is then calculated as:

$$\underline{t}_i^{dep}(k) = \begin{cases} \underline{t}_i^{p2}(k),\ if\ i \in \underline{I}_{bl}(k) \cap \mathcal{I}_{bl}^{real} \\ \phi_{r_i}(k) \cdot \underline{t}_i^{p2}(k) + \left(1 - \phi_{r_i}(k)\right) \cdot \underline{t}_i^{pre}(k), otherwise \end{cases}, \forall i \in \underline{I}_{bl}(k), k \in K \tag{32}$$

The speed of vehicle $i$ at $x_w$, $\underline{v}_i^{dep}(k)$, is determined as:

$$\underline{v}_i^{p1}(k) = v_i + a_i^{max} \cdot (\underline{t}_i^{dep}(k) - k), \forall i \in \underline{I}_{bl}(k), k \in K \tag{33}$$

$$\underline{v}_i^{p2}(k) = \underline{v}_i^{pre}(k) + a_i^{max} \cdot \hat{\tau}_i, \forall i \in \underline{I}_{bl}(k), k \in K \tag{34}$$

$$\underline{v}_i^{p3}(k) = min( v_i^{max}, \underline{v}_i^{p1}(k), \underline{v}_i^{p2}(k)), \forall i \in \underline{I}_{bl}(k), k \in K \tag{35}$$

$$\underline{v}_i^{dep}(k) = \begin{cases} \underline{v}_i^{p3}(k), if\ i \in \underline{I}_{bl}(k) \cap \mathcal{I}_{bl}^{real} \\ \phi_{r_i}(k) \cdot \underline{v}_i^{p3}(k) + (1 - \phi_{r_i}(k)) \cdot \underline{v}_i^{pre}(k), if\ i \in \underline{I}_{bl}(k) \cap \mathcal{I}_{bl}^{vir} \end{cases}, k \in K \tag{36}$$

Subsequently, we can proceed with the elaboration of $t_i^{dep}(k)$ by utilizing both $\underline{t}_i^{dep}(k)$ and $\underline{v}_i^{dep}(k)$ in further calculations. Using both $\underline{t}_i^{dep}(k)$ and $\underline{v}_i^{dep}(k)$, we can now calculate $t_i^{dep}(k)$. The earliest possible time $t_i^{p1}(k)$ for vehicle $i$ to pass the stop bar is:

$$t_i^{p1}(k) = \underline{t}_i^{dep}(k) + Z(x_c - x_w; \underline{v}_i^{dep}(k); v_i^{max}; a_i^{max}), \forall i \in \underline{I}_{bl}^{T}(k), k \in K \tag{37}$$

where $\underline{I}_{bl}^{T}(k)$ is a subset of $\underline{I}_{bl}(k)$ consisting of through vehicles.

The following formulas define the components of $Z(x_c - x_w; v_i(k); v_i^{max}; a_i^{max})$:

$$t_i^{acc}(k) = \frac{v_i^{max} - \underline{v}_i^{dep}(k)}{a_i^{max}}, \forall i \in \underline{I}_{bl}^{T}(k), k \in K \tag{38}$$

$$s_i^{acc}(k) = \frac{v_i^{max} + \underline{v}_i^{dep}(k)}{2} \cdot t_i^{acc}(k), \forall i \in \underline{I}_{bl}^{T}(k), k \in K \tag{39}$$

$$t_i^{p1}(k) = \begin{cases} \dfrac{-\underline{v}_i^{dep}(k) + \sqrt{(\underline{v}_i^{dep}(k))^2 + 2a_i^{max}(x_c - x_w)}}{a_i^{max}},\ if\ s_i^{acc}(k) > x_c - x_w \\ \dfrac{x_c - x_w - s_i^{acc}(k)}{v_i^{max}} + t_i^{acc}(k), otherwise \end{cases}, \forall i \in \underline{I}_{bl}^{T}(k), k \in K \tag{40}$$

If a right-turning vehicle changes lanes, the car-following relationships are reconfiguration. Let $\bar{i}'$ denote the nearest through vehicle ahead of $i$. The passing time $t_i^{p2}(k)$ is then:

$$t_i^{p2}(k) = max( t_i^{p1}(k), t_{\bar{i}'}^{dep}(k) + \hat{\tau}_i), \forall i \in \underline{I}_{bl}^{T}(k), k \in K \tag{41}$$

The constraints outlined in Eqs.(27)-(29) also apply to vehicles in $\underline{I}_{bl}^{T}(k)$ for deriving $t_i^{dep}(k)$.

### *3.2.2 Estimation of passing states for vehicle in general lane*

(1) Estimation of the departure time $t_i^{dep}(k)$ for vehicles in set $I'_{gl}(k)$:

For vehicle $i$ that has already passed the stop bar at time $k$ in the general lane, the $t_i^{dep}(k)$ can directly be retrieved from the known vehicle state set $S$, i.e., $t_i^{dep}(k) = t_i^{tra}(k)$.

(2) Estimation of the departure time $t_i^{dep}(k)$ for vehicles in set $\overline{I}_{gl}(k)$:

For vehicles that have not yet passed the stop bar at time $k$, the constraints specified in Eqs.(22) -(29) also apply. However, unlike CAVs, HDVs experience additional start-up loss time when they arrive at the stop bar during a red phase and begin moving at the onset of the green phase. Consequently, the passing time $t_i^{p3}(k)$ for these vehicles is calculated as:

$$t_i^{p3}(k) = max(t_i^{p2}(k), t_R + t_r + \tau_r + \tau_a), \forall i \in \overline{I}_{gl}(k) \cap \mathcal{I}_h, k \in K \tag{42}$$

where $\tau_r$ is the reaction time to the green light, $\tau_a$ is the time required to start and accelerate to pass the stop bar.

Specifically, considering that CAVs may change lanes into bus lane, the passing time $t_i^{dep}(k)$ for vehicles in set $\overline{I}_{gl}(t) \cap \mathcal{I}_a$ is calculated as:

$$t_i^{dep}(k) = \left(1 - \phi_i(k)\right) \cdot t_i^{p3}(k) + \phi_i(k) \cdot t_i^{pre}(k) + \left(1 - \lambda(k)\right) \cdot M, \forall i \in \overline{I}_{gl}(t) \cap \mathcal{I}_a, k \in K \tag{43}$$

### *3.3 Modelling the minimum time headway of different types of vehicles*

In the proposed model, Newell's car-following model is used to characterize the minimum time headway of CAVs, HDVs, and CABs when closely following the preceding vehicle. Specifically, in Newell's framework, the distance between the current position of vehicle $i$ and the position of its preceding vehicle $\overline{i}$ after a time delay $\tau_i$ is represented by a constant spacing $d_i$. This relationship is expressed as:

$$x_i(k + \tau_i) = x_{\overline{i}}(k) - d_i \tag{44}$$

where $\tau_i$ and $d_i$ represent the temporal and spatial displacements in Newell's car-following model.

It's important to note that many studies employing Newell's model assume uniform vehicle dimensions. However, our study faces a unique challenge due to the presence of CABs, which vary significantly in size compared to private cars. As shown in Fig. 9, based on Newell's car-following model (Wei et al. 2017), $d_i$ is influenced by both the type of the preceding vehicle and the ego vehicle. To better reflect the distinctions, we have refined $\tau_i$ and $d_i$ specifically for HDVs, CAVs, and CABs. The value of $\tau_i$ depends on the type of vehicle $i$, while $d_i$ depends on both the type of vehicle $i$ and its preceding vehicle $\overline{i}$.

$$\tau_i = \begin{cases} \tau_A, when\ i \in \mathcal{I}_b \cup \mathcal{I}_a \\ \tau_H, when\ i \in \mathcal{I}_h \end{cases} \tag{45}$$

$$d_i = d'_i + l_{\overline{i}} \tag{46}$$

$$d'_i = \begin{cases} d_A, when\ i \in \mathcal{I}_b \cup \mathcal{I}_a \\ d_H, when\ i \in \mathcal{I}_h \end{cases} \tag{47}$$

$$l_{\bar{i}}, = \begin{cases} l_B, when\ \bar{i} \in \mathcal{I}_b \\ l_V, when\ \bar{i} \in \mathcal{I}_a \cup \mathcal{I}_h \end{cases} \tag{48}$$

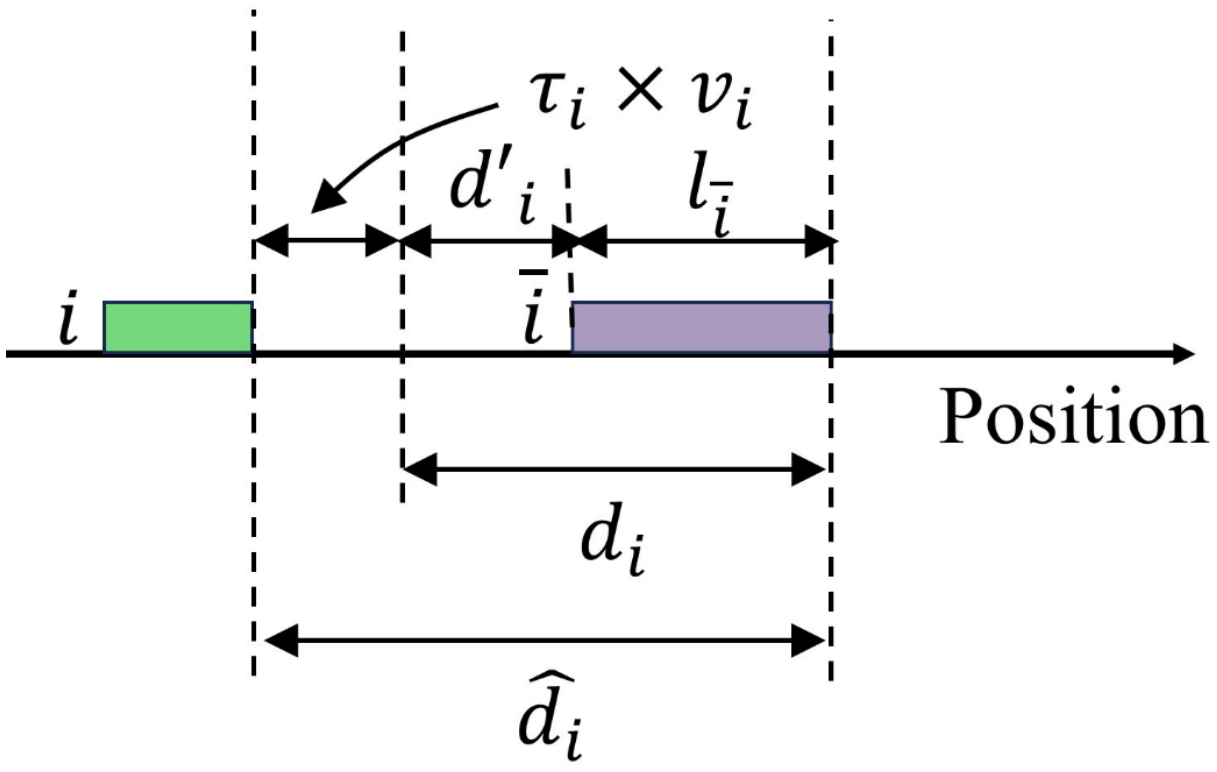

Fig. 9 The interconvertible calculation between time headway and space headway

Furthermore, when vehicle $i$ is closely following vehicle $\bar{i}$, the minimum time headway $\hat{\tau}_i$ and minimum safety distance $\hat{d}_i$ can be calculated as follows:

$$\hat{\tau}_i = \tau_i + \frac{d_i}{v_i(k)} \tag{49}$$

$$\hat{d}_i = \tau_i \cdot v_i(k) + d_i \tag{50}$$

We assume that each vehicle aims to pass the stop bar as quickly as possible. Therefore, a vehicle will attempt to accelerate to its maximum speed and maintain that speed until it is constrained by the preceding vehicle, entering a state of close following. In this case, the time headway satisfies Eq.(49) (Dai et al. 2023).

### *3.4 Solution method*

#### *3.4.1 Defining the Extent of ROW optimization*

When a bus is dwelling at a bus stop, granting access to private cars trailing behind it into the bus lane may cause additional parking delays for those vehicles and hinder the timely arrival and parking of subsequent buses. Moreover, the unpredictability of bus dwell time, which is inherently tied to passenger volume, makes it difficult to accurately estimate bus travel times during dwelling. In such cases, if too many private cars are allowed to enter the bus lane ahead of the bus just before it is about to depart, it may excessively delay the bus from passing the stop bar.

Therefore, before applying the optimization model, we carefully define the sets $\mathcal{I}_{gl}$ and $\mathcal{I}_{bl}$ to be optimized, ensuring that the aforementioned issues are avoided:

**Step 1:** Initialize the current time $k_0$ and update the states of all vehicles in both the general and bus lanes. Define $\mathcal{I}_{gl}$ as the set of vehicles in the general lane that are located between the stop bar and the nearest right-turning vehicle (with the position of this vehicle defined as $x_r(k_0)$). Define $\mathcal{I}_{bl}$ as the set of all vehicles in the bus lane.

**Step 2:** If a bus $j$ is detected in a dwelling state, calculate the earliest time for the bus to start and reach the stop bar as $t_1$=$k_0 + Z\left(x_c - x_s; v_j(k_0); v_j^{max}; a_j^{max}\right)$. Update the sets:

$$\mathcal{I}_{gl} = \left\{i | x_s \le x_i(k_0) \le x_c, i \in \mathcal{I}_{gl}\right\} \tag{51}$$

$$\mathcal{I}_{bl} = \{i | \max(x_s, x_r(k_0)) \le x_i(k_0) \le x_c, i \in \mathcal{I}_{bl}\} \tag{52}$$

Here, $x_s$ denotes the position of bus stop. Proceed to Step 4. If no bus is detected, continue to Step 3.

**Step 3:** If no dwelling bus is detected, set $t_1 = M$ (a large number). If there are buses in the segment $[0, x_s]$, identify the bus closest to the stop as index $j$ and update:

$$\mathcal{I}_{bl} = \left\{i | x_j(k_0) \le x_i(k_0) \le x_c, i \in \mathcal{I}_{bl}\right\} \tag{53}$$

$$\mathcal{I}_{gl} = \left\{i | \max(x_j(k_0), x_r(k_0)) \le x_i(k_0) \le x_c, i \in \mathcal{I}_{gl}\right\} \tag{54}$$

Proceed to Step 4.

**Step 4:** Iterate through each vehicle in $\mathcal{I}_{gl}$. It is assumed that the bus will start at $k_0$. If the $t_i^{p1}(k_0)$ of a vehicle $i$, plus the minimum time headway of the bus $\hat{\tau}_b$, does not exceed the time it takes for the bus to reach the stop bar, i.e.,

$$t_i^{p1}(k_0) + \hat{\tau}_b > t_1 \tag{55}$$

Update $\mathcal{I}_{gl} \leftarrow \mathcal{I}_{gl} \backslash \{i\}$. Otherwise, retain the current set $\mathcal{I}_{gl}$.

**Step 5:** Output the updated sets $\mathcal{I}_{gl}$ and $\mathcal{I}_{bl}$.

### *3.4.2 ROW Pre-allocation Heuristic Algorithm*

This section introduces the ROW Pre-allocation Heuristic (ROWPH) algorithm, which applies heuristic principles to reduce problem complexity. The $\vartheta_i(k)$ acts as an indicator of lane-changing opportunities for vehicles in the general lane, helping streamline the problem by reducing constraints. Specifically, $\vartheta_i(k) = 1$ signifies a viable lane-changing opportunity for vehicle $i$ at time $k$, while $\vartheta_i(k) = 0$ means no such opportunity exists. Additionally, the algorithm employs dimensionality reduction techniques to enhance computational efficiency by eliminating redundant dimensions.

**Step 1:** Initialize the start time $k = k_0$ and update the vehicle states $S$. Set $\vartheta_i(k) = 0$ for $i \in \mathcal{I}_{gl} \cap \mathcal{I}_a$, $k \in K$.

**Step 2:** For each vehicle $i$ in the set $\mathcal{I}_{bl}^{real}\backslash I'_{bl}(k)$, compute: $x_{lead}(k) = x_{\bar{i}}(k) - l_{\bar{i}}$ , $x_{follow}(k) = x_i(k)$ , $d_{gap}(k) = x_{lead}(k) - x_{follow}(k)$ . If $d_{gap}(k) > 2d_{safe} + l_v$, add the interval $[x_{follow}, x_{lead}]$ to the list of available spaces.

- If $i$ is the first vehicle in the bus lane at time $k$, set $x_{lead} = x_c$.
- If $i$ is the last vehicle, account for both the space ahead and the space between the vehicle and the road entrance, setting $x_{lead} = x_i(k) - l_i$ and $x_{follow} = 0$.

**Step 3:** Check available spaces. For each space $[x_{follow}, x_{lead}]$, evaluate vehicles $i \in \mathcal{I}_{gl} \cap \mathcal{I}_a$. If conditions in Eqs.(56)-(59) are satisfied, set $\vartheta_i(k) = 1$; otherwise, $\vartheta_i(k) = 0$.

$$x_i(k) < x_n \tag{56}$$

$$v_i(k) > 0 \tag{57}$$

$$x_{follow} + d_{safe} + l_i < x_i(k) < x_{lead} - d_{safe} \tag{58}$$

$$x_{\underline{r}_i}(k + k_{lc}) < x_i(k + k_{lc}) < x_{\bar{r}_i}(k + k_{lc}) \tag{59}$$

Consequently, constraints Eq.(10)-(15) can be superseded by constraint Eq. (60).

$$\phi_i(k) \le \vartheta_i(k), \forall i \in \mathcal{I}_{gl} \cap \mathcal{I}_a, k \in K \tag{60}$$

Proceed to Step 4.

**Step 4:** If $k > k_0$ and the vehicles with $\vartheta_i(k - \Delta k) = 1$ remain the same as those with $\vartheta_i(k) = 1$, remove $k$ from the set $K$ (i.e., $K \leftarrow K\backslash\{k\}$). Proceed to Step 5.

**Step 5:** Increment $k$ by $\Delta k$. If $k \in K$, repeat Step 2. Otherwise, output $\vartheta_i(k)$ and the updated $K$.

### *3.4.3 Model Linearization*

To transform nonlinear constraints into linear forms, we linearize the model as a MILP problem. Specifically, the nonlinear constraint presented in Eq. (28) can be reformulated as:

$$t_i^{flo}(k) \le \frac{t_i^{p2}(k)}{t_c} \tag{61}$$

$$t_i^{flo}(k) \ge \frac{t_i^{p2}(k)}{t_c} - 1 + P \tag{62}$$

$$t_i^R(k) = t_i^{flo}(k) \cdot t_c \tag{63}$$

where $t_i^{flo}(k)$ is an auxiliary variable indicating the signal cycle during which vehicle $i$ arrives at the stop bar. When $\frac{t_i^{p2}(k)}{t_c}$ is non-integer, $t_i^{flo}(k)$ takes the largest integer less than this value, leveraging integer programming properties. $P$ is a very small positive number.

To address the nonlinearity in Eq. (40), we introduce the binary auxiliary variable $\sigma_i(k)$, as defined in Eqs.(64)-(67). This variable serves as an indicator:

- $\sigma_i(k) = 1$: When $x_c - x_w \le s_i^{acc}(k)$, meaning the travel time of vehicle $i$ consists solely of the acceleration phase.
- $\sigma_i(k) = 0$ when $x_c - x_w \ge s_i^{acc}(k)$, indicating that the travel time includes both acceleration and uniform cruising phases.

The following equations define the conditions for $\sigma_i(k)$:

$$x_c - x_w - s_i^{acc}(k) \le M \cdot (1 - \sigma_i(k)) \tag{64}$$

$$x_c - x_w - \underline{v}_i^{dep}(k) \cdot (t_i^{p1}(k) - k) - \frac{a_i^{max} \cdot (t_i^{p1}(k) - k)^2}{2} = M \cdot (1 - \sigma_i(k)) \tag{65}$$

$$s_i^{acc}(k) - (x_c - x_w) \le M \cdot \sigma_i(k) \tag{66}$$

$$t_i^{p1}(k) - \left( \frac{(x_c - x_w) - s_i^{acc}(k)}{v_i^{max}} - t_i^{acc}(k) \right) = M \cdot \sigma_i(k) \tag{67}$$

The linearization ensures that the entire model becomes a MILP problem, which is solved using Gurobi on a computer with an Intel Core i7-12700 CPU, 16GB RAM, running Windows 10 64-bit.

## *3.5 Dynamic control scheme*

Due to the inherent unpredictability in HDV movements, using driving models to estimate vehicles' passing states within the optimization framework can lead to inaccuracies. As a result, some CAVs granted the ROW may not have the appropriate conditions to execute a lane change in practice. To manage these uncertainties and improve system resilience against modeling errors, we introduce a rolling horizon approach combined with an activation mechanism, known as the dynamic adaptive ROW allocation protocol. This protocol dynamically implements the ROW allocation optimization model, allowing for real-time adaptability to changing traffic conditions. Table 1 presents the pseudo-code for the protocol, $k_c$ denotes the lane change time communicated by the control center to vehicles following the ROW optimization. $\mathcal{I}_c$ represents the set of vehicles granted ROW, i.e., $\mathcal{I}_c = \{i | \phi_i(k_c) = 1, i \in \mathcal{I}_{gl}\}$. $p_i$ refers to the anticipated preceding vehicle in the bus lane for vehicle $i$ that has been granted ROW.

Table 1 Dynamic adaptive ROW allocation protocol

| **Input:** Current time $k_0$, vehicle states $S$, set of vehicles granted ROW $\mathcal{I}_c$, lane change time $k_c$. |
|---|
| At any given current time $k_0$, the control center needs to perform the following operations |
|  **If** $\mathcal{I}_c$ is empty: |
|   Execute the ROW allocation optimization model |
|   Update $\mathcal{I}_c$ and $k_c$ |
|   Send lane change recommendations to vehicles in set $\mathcal{I}_c$ |
|  **Else**: |
|   **If** $k_0 > k_c$: |
|    **For** each vehicle $i$ in $\mathcal{I}_c$: |
|     **If** $i$ is already in the bus lane: |

Remove $i$ from $\mathcal{I}_c$
**Else**:
Calculate $d_p$, $d_p = x_{\overline{r_i}}(k_0) - x_i(k_0) - d_{safe}$
Calculate $d_f$, $d_f = x_i(k_0) - x_{\underline{r_i}}(k_0) - d_{safe}$
**If** unsafe conditions exist ($d_p < 0$ or $d_f < 0$) or the intended preceding vehicle $p_i$ does not match the actual preceding vehicle $\overline{r_i}$:
Empty $\mathcal{I}_c$ and cancel the remaining ROW for vehicles that have not yet entered the bus lane
**Else**:
Retain the current $\mathcal{I}_c$
**End**
**Else**:
Retain the current $\mathcal{I}_c$

Furthermore, a system diagram has been developed to illustrate the driving process for all vehicles traversing a road segment featuring a bus stop and a signalized intersection, as depicted in Fig. 10.

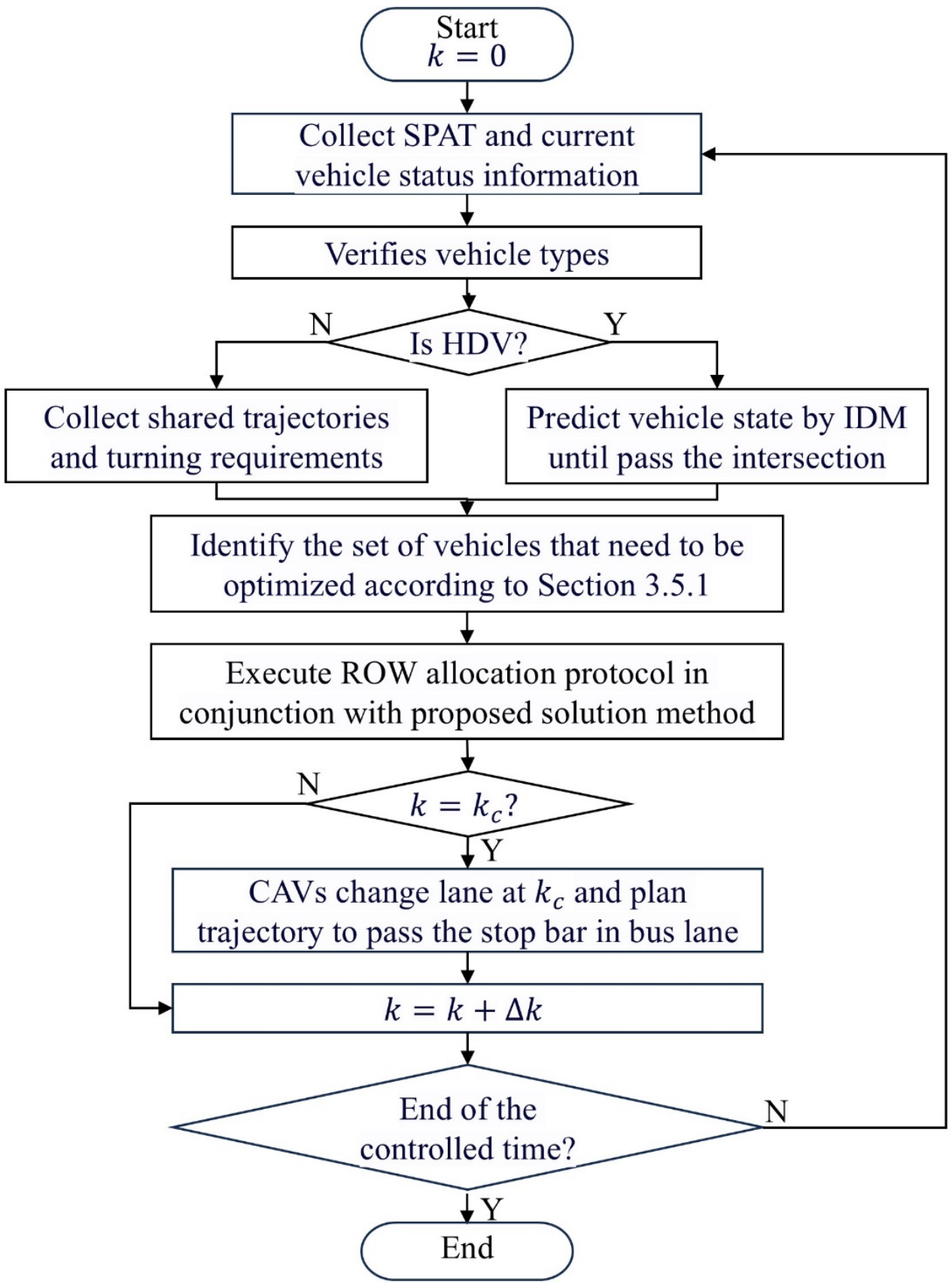

Fig. 10 Dynamic control framework

## 4. Numerical examples

### *4.1 Experimental setups*

In this section, we conduct simulation experiments using SUMO, focusing on two distinct scenarios to evaluate the effectiveness of the proposed strategy.

- Scenario A: Without a right-turn pocket. This scenario assesses how the proposed strategy

utilizes bus lanes to ease congestion in general lanes.

- Scenario B: With a right-turn pocket. This scenario examines the applicability and effectiveness of the strategy in a complex environment where both general and bus lanes accommodate a mix of right-turning and through vehicles. This scenario also explores the temporary use of the bus lane by right-turning vehicles before they transition into the right-turn pocket.

Our longitudinal trajectory planning model for CAVs follows the method proposed by Ma et al. (2021). Unlike CAVs, CABs execute their trajectory planning in three distinct phases: (1) approaching the bus stop, (2) dwelling at the stop, and (3) passing the stop bar. The method also considers the combined influence of bus stop locations and traffic signals, inspired by the work of Shan et al. (2024). Additionally, to demonstrate adaptability, we use the stochastic Krauss model in SUMO to simulate HDV driving behavior, although our method employs the IDM model for HDV trajectory prediction.

The parameter settings used in the simulations are summarized in Table 2. Notably, the bus stop is designed to accommodate up to two buses dwelling simultaneously. Five random seeds are used to ensure robustness, with each simulation running for 1800 seconds, including a full warm-up period (Li et al. 2023). The resolution of the simulator is 1 s.

Our analysis begins with benchmark scenarios: (1) the private car demand $Q_{veh}$ is 720 vehicles per hour (veh/h); (2) bus arrival intervals ($F$) are 60 seconds, with a standard deviation of 20 seconds; (3) the distance between the road entry and the bus stop ($x_s$) is 150 meters; and (4) the right-turn ratio in Scenario B is 0.2. Additionally, we perform a sensitivity analysis to assess how key factors—such as $Q_{veh}$, $F$, $x_s$and right-turn ratio—influence the performance of the DBPL control strategy.

Table 2 Simulation parameters setup

| Parameter | Value |
|---|---|
| CAB length $l_B$ (m) | 8 |
| CAV/HDV length $l_V$ (m) | 4 |
| Maximum vehicle speed $v_i^{max}$ (m/s) | 14 |
| Max/Min acceleration $a_i^{max/min}$ (m/s$^2$) | 2 |
| Perception and reaction time of CAV/CAB $\tau_A$ (s) | 1 |
| Perception and reaction time of HDV $\tau_H$ (s) | 2 |
| Safety buffer (CAV/CAB) $d_A$ (m) | 1.5 |
| Safety buffer (HDV) $d_H$ (m) | 2.5 |
| Time to start and pass stop bar $\tau_a$ (s) | 1.5 |
| Green light reaction time $\tau_r$ (s) | 0.4 |
| Lane-change safe distance $d_{safe}$ (m) | 6 |
| **Road parameters** | |
| No-changing zone length (m) | 30 |
| Distance from entry to bus stop $x_s$ (m) | 150 |
| Control zone length (m) | 400 |
| Through car demand $Q_{veh}$ (veh/h) | 720 |
| Bus arrival interval mean (s) | 60 |
| Bus arrival interval variance (s) | 20 |
| Bus dwelling time mean (s) | 30 |
| Bus dwelling time variance (s) | 20 |
| Right-turn pocket length (scenario B) (m) | 130 |
| **Signal phase information** | |
| Signal cycle length $t_c$(s) | 60 |
| Green phase length $t_b$ (s) | 30 |
| Red and amber phase length $t_r$ (s) | 30 |

### *4.2 benchmark scenario analysis*

This section evaluates the performance of two strategies, EBL and DBPL, across various CAV Market Penetration Rates (MPRs).

Fig. 11 shows the average travel time of private cars in both scenarios, along with the reduction achieved by DBPL compared to EBL. As the CAV MPR increases, both strategies yield shorter travel times due to CAVs' efficient intersection traversal. However, DBPL achieves a more substantial reduction, demonstrating the benefit of allowing certain CAVs to use bus lanes to alleviate general lane congestion. Additionally, it can be observed that the reduction in travel time initially increases and then decreases, indicating that the DBPL strategy is most effective in improving mixed traffic conditions at moderate CAV MPRs. For instance, at an MPR of 40%, travel times in Scenario A and Scenario B are reduced by approximately 20% and 8%, respectively. The reason for this trend is that, as the number of CAVs on the road increases, the overall traffic efficiency improves even without implementing the bus lane borrowing strategy. However, as MPRs continue to rise and traffic pressure diminishes, the DBPL strategy still provides some degree of additional improvement.

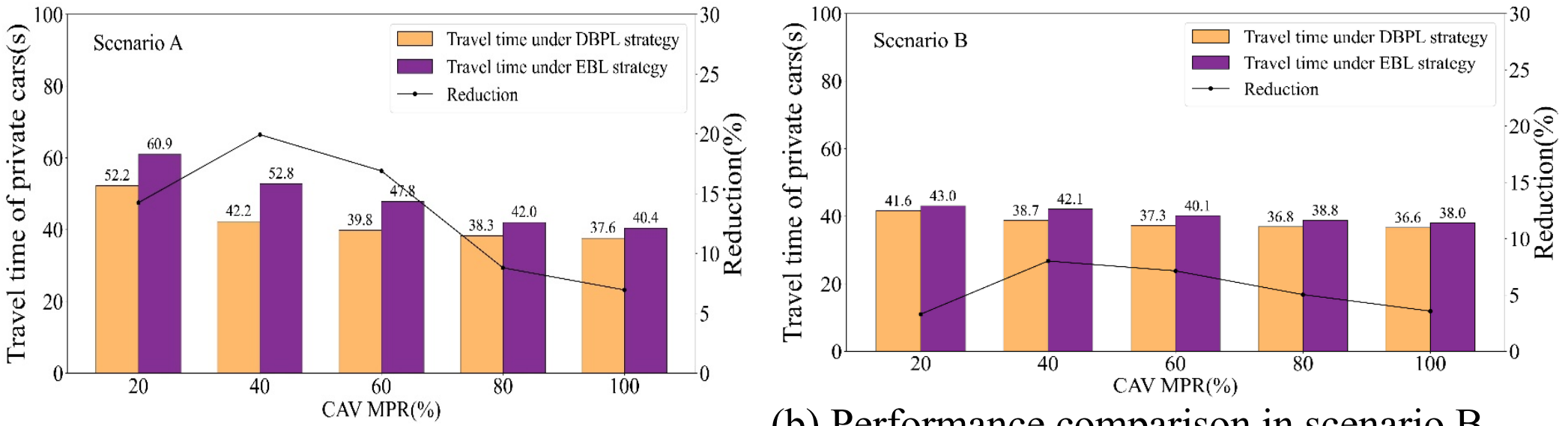

(a) Performance comparison in scenario A

(b) Performance comparison in scenario B

Fig. 11 Average travel time and reduction of private car in scenario A and B

Fig. 12 (a) illustrates that DBPL outperforms EBL in reducing travel times across all vehicle types in Scenario A. In this figure, A, H, and B represent CAVs, HDVs, and CABs, respectively.

CAB travel times remain almost identical between the two strategies, with a maximum difference of only 0.3 seconds, ensuring bus priority is fully maintained. Private cars benefit significantly. At 40% MPR, CAV travel times decrease by 13 seconds and HDV travel times by 9 seconds under DBPL, compared to EBL. This improvement stems from two factors: CAVs can pass through the stop bar faster, and general lane congestion is eased, indirectly benefiting HDVs. Fig. 12 (b) shows that in Scenario B, DBPL maintains its advantage despite the presence of right-turn traffic. Notably, at 40% MPR, CAV travel times are reduced by around 5 seconds, and HDV times by 3 seconds. Additionally, the average travel time of private cars in Scenario B is generally shorter than in Scenario A. This is because right-turning vehicles in Scenario B experience shorter travel times compared to non-right-turning vehicles as shown in Fig. 13. Since right-turning vehicles are not controlled by traffic signals and have temporary access to the bus lane, they can directly enter the right-turn pocket and pass through the intersection smoothly.

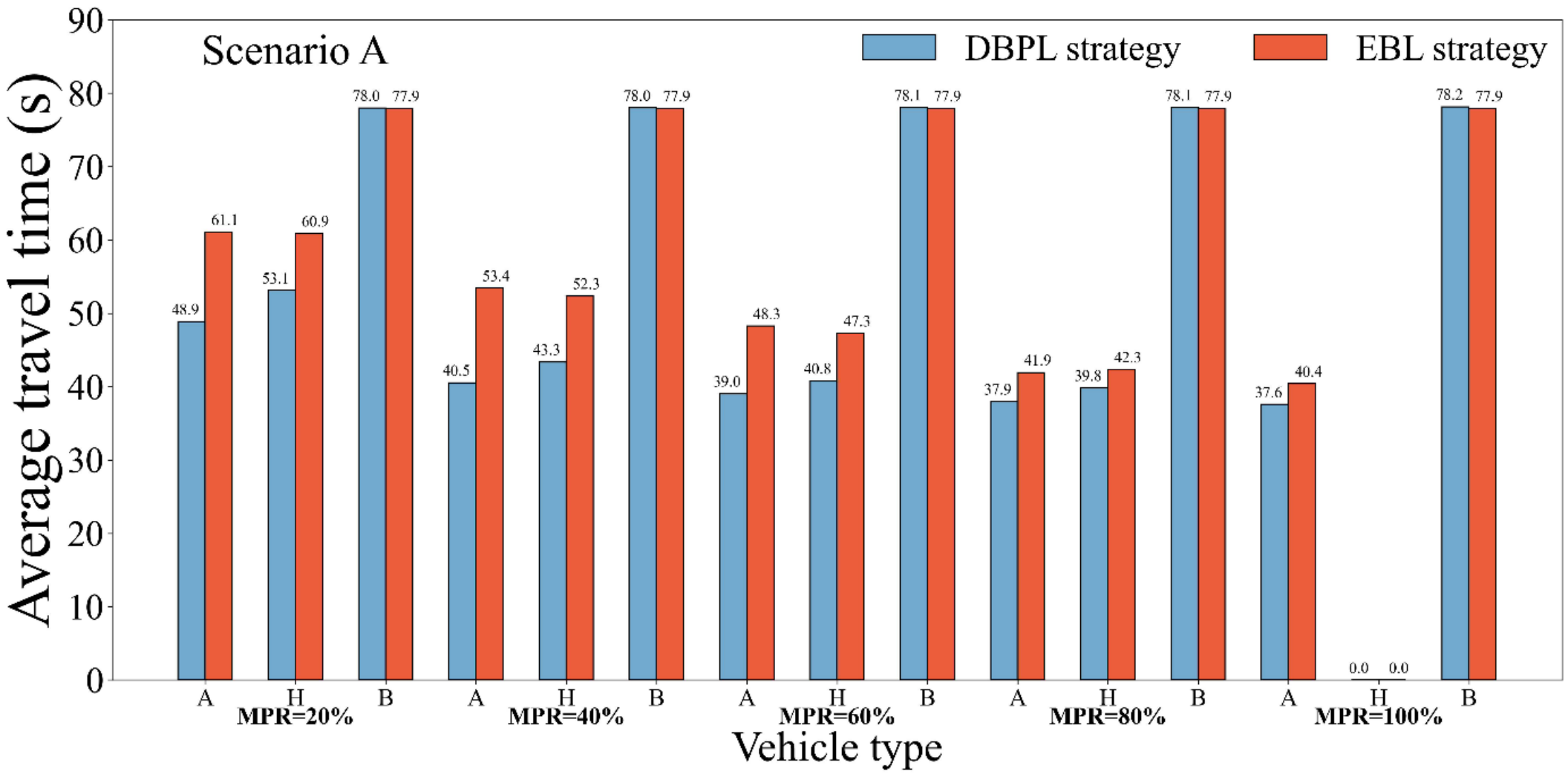

(a) Scenario A

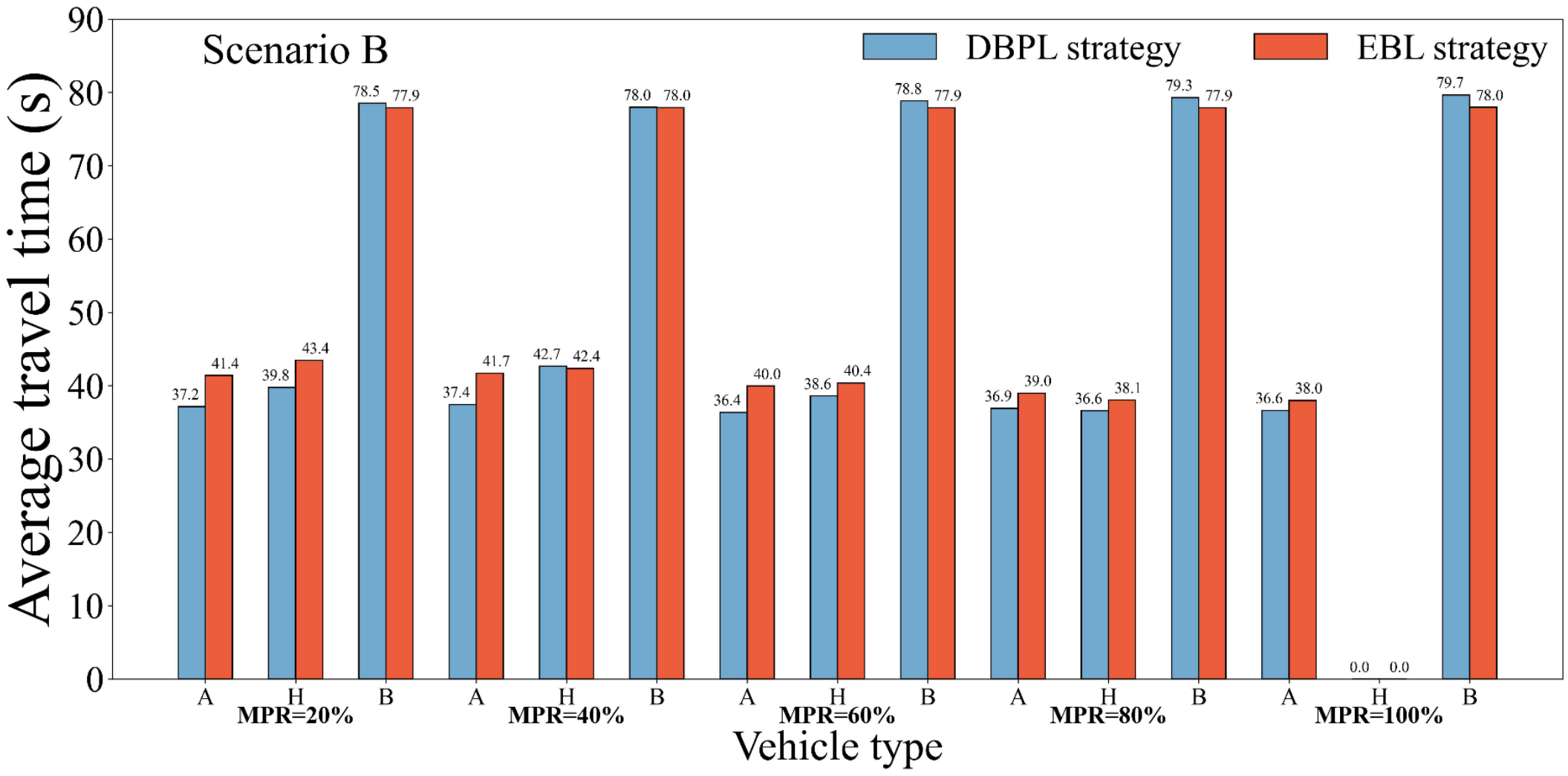

(b) Scenario B

Fig. 12 Average travel time grouped by vehicle type under different CAV MPRs

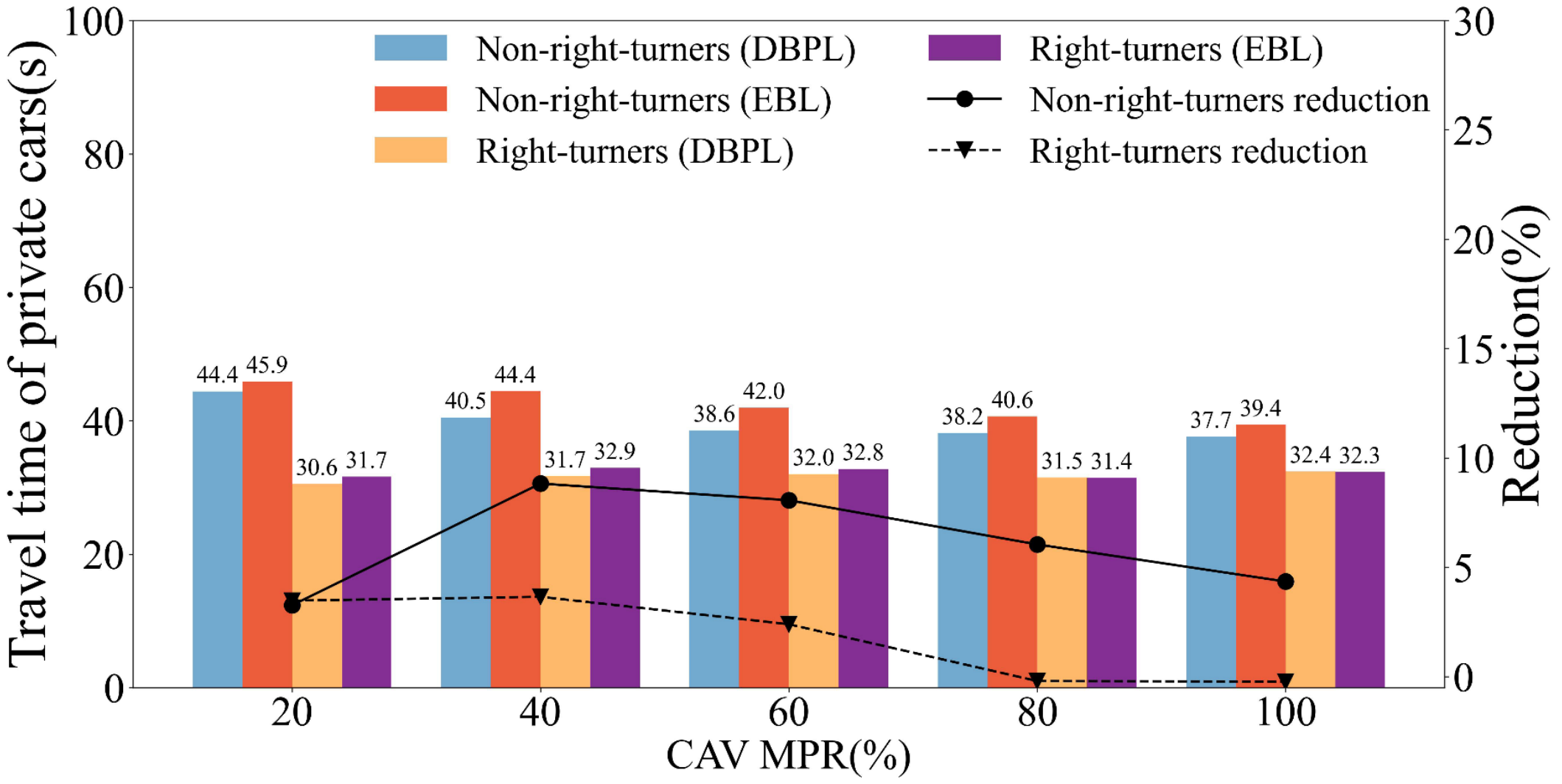

Fig. 13 Average travel time grouped by vehicle movements under different CAV MPRs in scenario B

Fig. 14 compares the spatio-temporal trajectories of vehicles under EBL and DBPL strategies in Scenario A at 40% MPR. The blue lines represent HDVs, while red lines depict CAVs' trajectories. As shown in Fig. 14 (a) and (c), under the EBL strategy, HDVs stop at the stop bar during the red phase, creating a shockwave that leads to start-up losses. In contrast, CAVs utilize SPaT information to pass through intersections smoothly, minimizing delays. Additionally, the leading effects of some CAVs become evident, helping to reduce the impact of the shockwave and mitigate start-up losses in mixed traffic. Moreover, Fig. 14 (c) reveals that there are sufficient available spatial-temporal resources for general traffic to utilize. In contrast, as shown in Fig. 14 (b), the DBPL strategy alleviates congestion in the general lane by redirecting some CAVs into the bus lane. However, not all CAVs switch lanes, allowing certain vehicles to remain in the general lane and preserve their leading effect, thereby smoothing traffic flows. Additionally, CAVs that enter the bus lane pass the stop bar earlier compared to those under the EBL strategy. A comparison between Fig. 14 (c) and Fig. 14 (d) further demonstrates that the proposed control strategy has minimal impact on bus priority, even with CAVs operating in the bus lane.

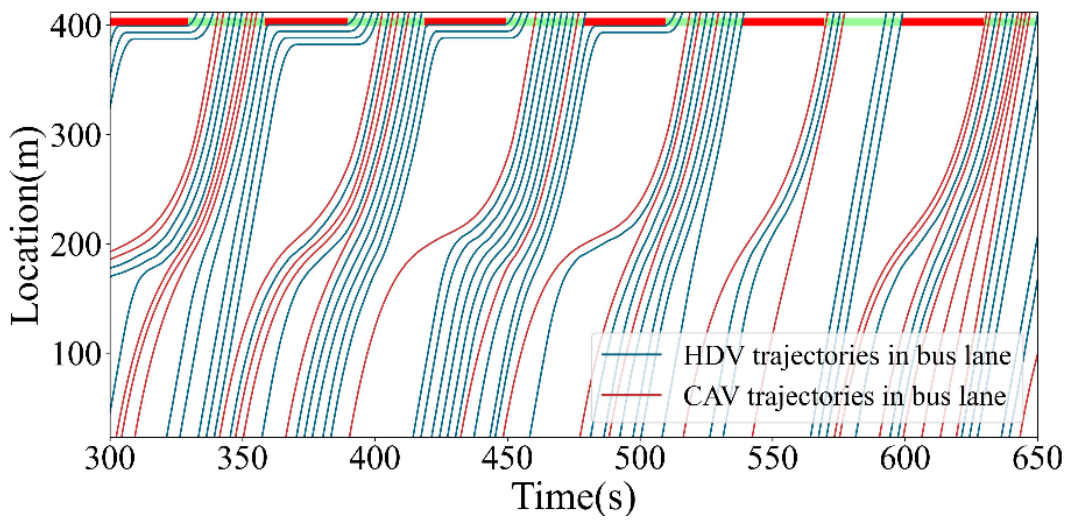

(a) Trajectory in general lane under EBL

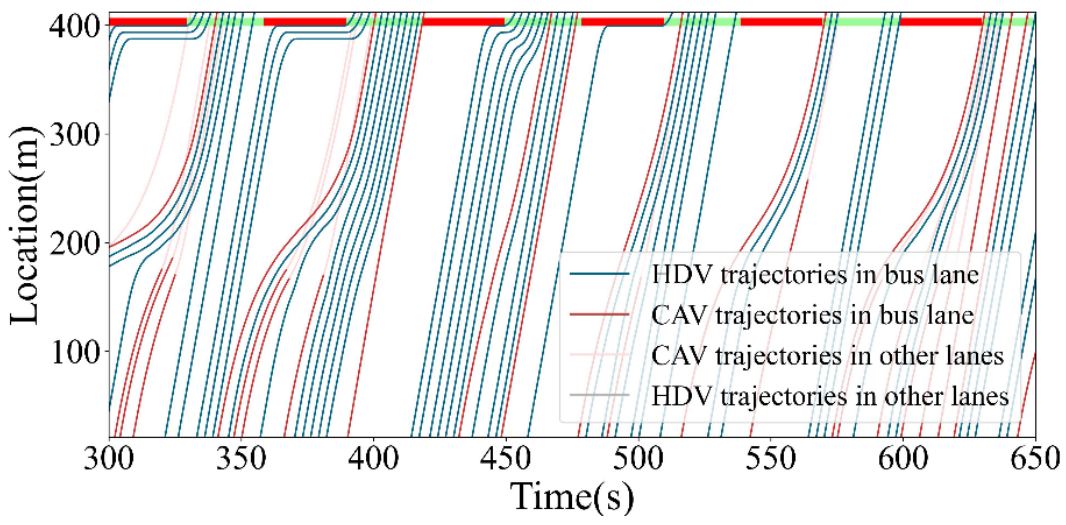

(b) Trajectory in general lane under DBPL

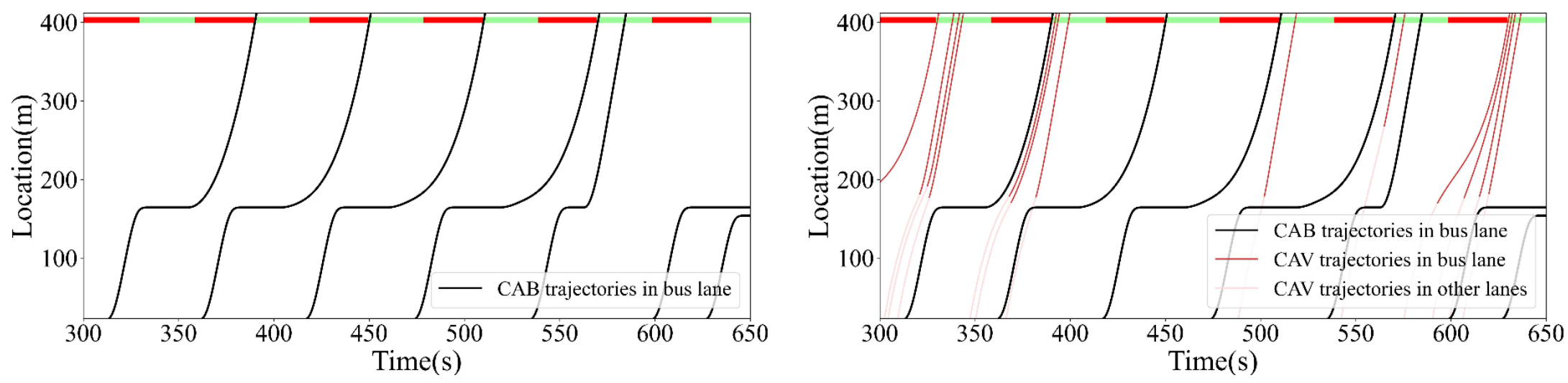

(c) Trajectory in bus lane under EBL (d) Trajectory in bus lane under DBPL

Fig. 14 Vehicle trajectories between DBPL and EBL in scenario A

Fig. 15 provides a detailed illustration of the benefits of the proposed strategy based on the vehicle trajectories in Scenario B (MPR=40%). It shows that right-turning HDVs and CAVs initially travel in the general lane, temporarily switch to the bus lane, and then enter the right-turn pocket to exit the intersection. When the demand for right-turn traffic is satisfied, some through CAVs also use the bus lane, significantly reducing travel times for vehicles in the general lane. For instance, comparing the third green phase in Fig. 15 (a) and Fig. 15 (b) reveals that after some CAVs enter the bus lane, the following HDVs can pass the stop bar more quickly. Moreover, these CAVs pass through the intersection earlier than they would if they stayed in the general lane. Similarly, a comparison between Fig. 15 (c) and Fig. 15 (d) shows that the DBPL strategy maintains bus trajectories almost identical to those under the EBL strategy, indicating minimal interference from private cars. Additionally, the lane-changing patterns of right-turning vehicles remain consistent, confirming that the DBPL strategy continues to ensure their smooth passage through the intersection.

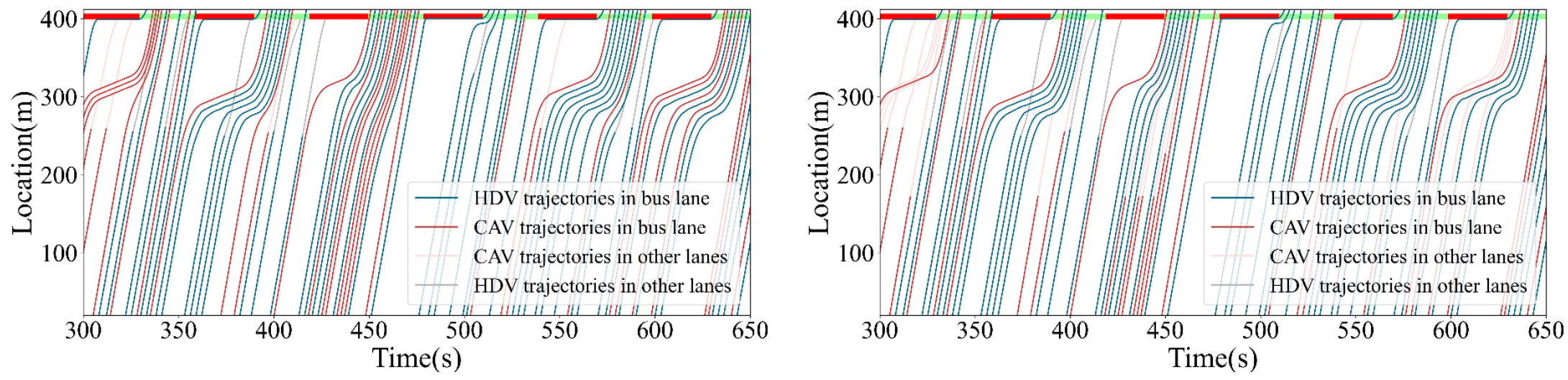

(a) Trajectory in general lane under EBL (b) Trajectory in general lane under DBPL

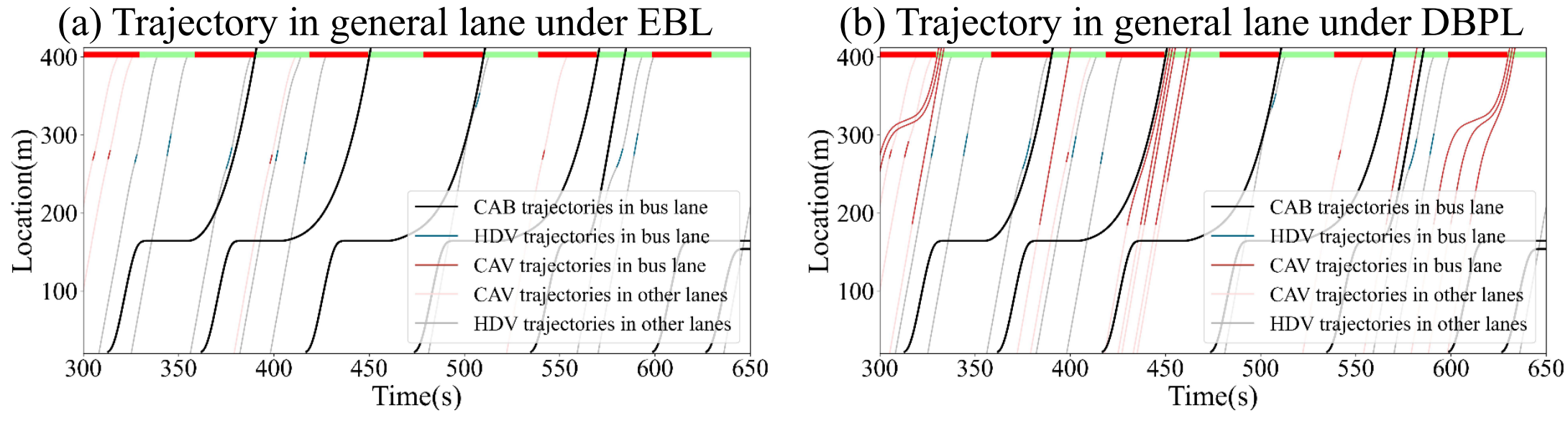

(c) Trajectory in bus lane under EBL (d) Trajectory in bus lane under DBPL

Fig. 15 Vehicle trajectories between DBPL and EBL in scenario B

### *4.3 Sensitivity analysis*

### *4.3.1 Sensitivity to private car demand*

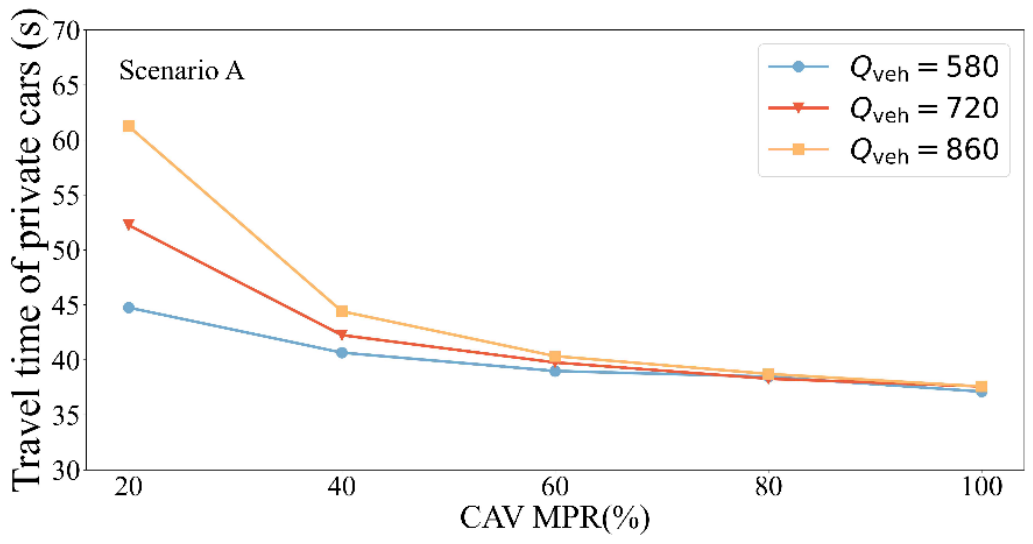

(a) Average travel time in scenario A

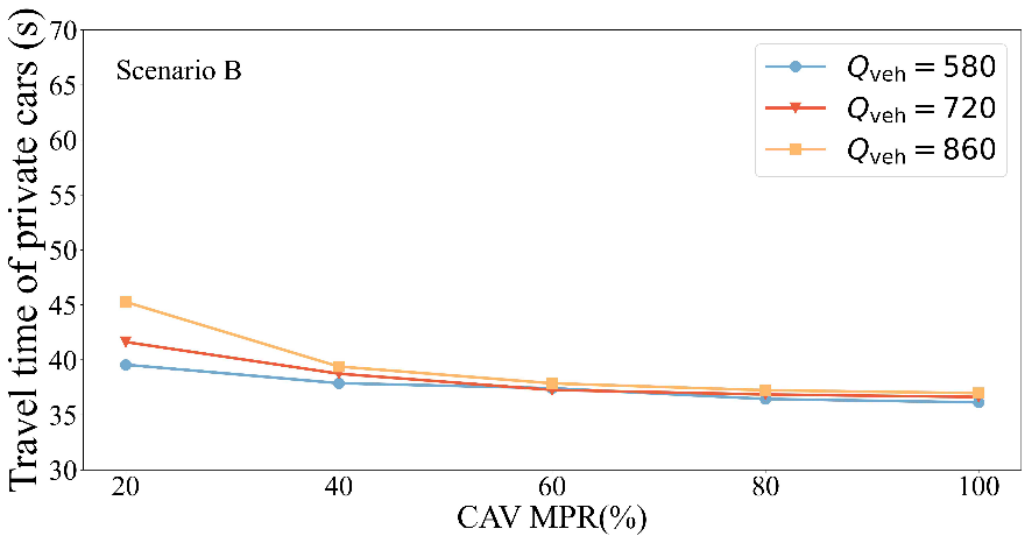

(b) Average travel time in scenario B

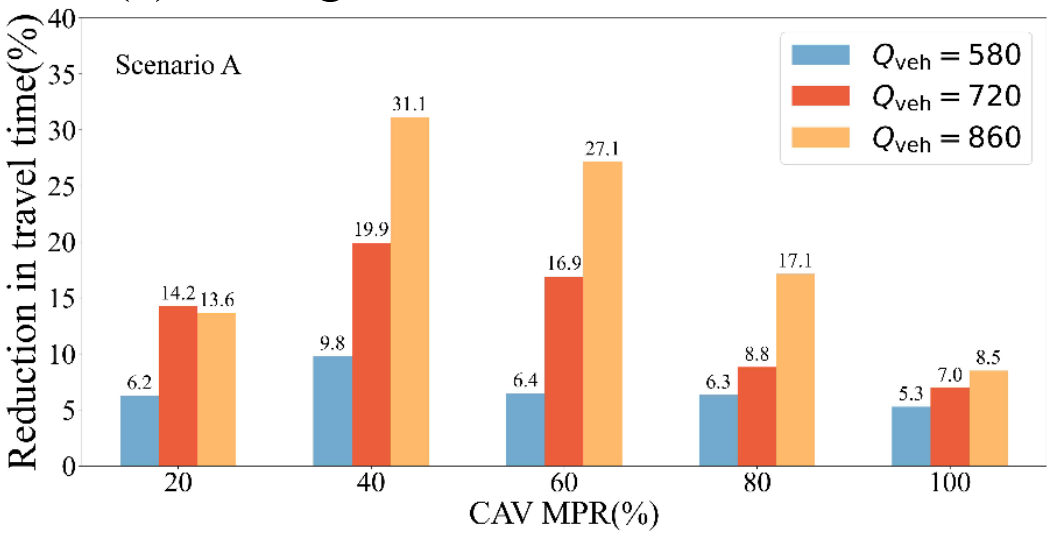

(c) Travel time reduction in scenario A

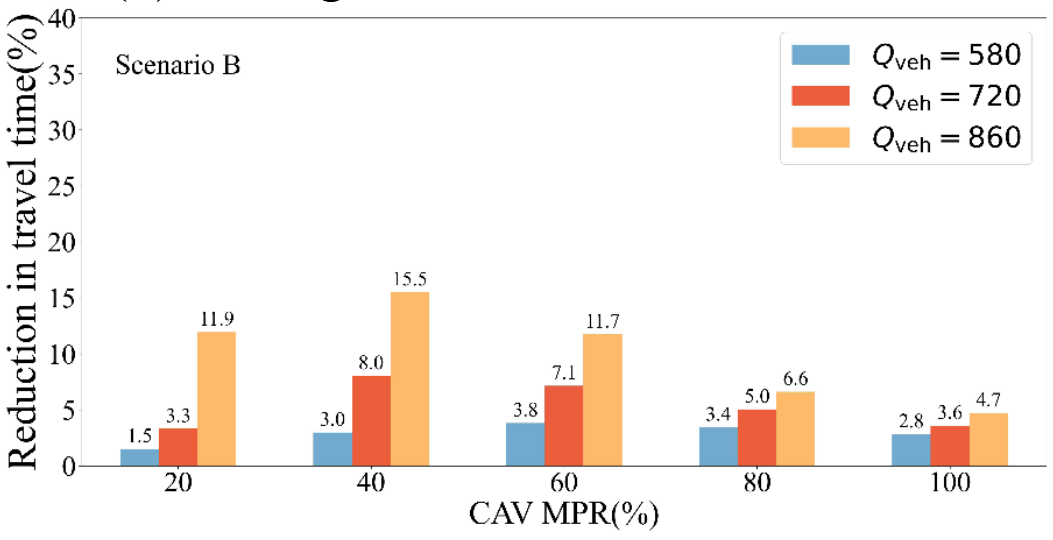

(d) Travel time reduction in scenario B

Fig. 16 Impacts on average travel time of the CAV MPRs at different demand levels.

Fig. 16 shows the impact of varying private car demand on average travel times under the DBPL strategy. General traffic volumes are set to 580, 720, and 860 veh/h to represent different demand levels. The results for Scenarios A and B are illustrated in Fig. 16 (a) and (b), respectively. The traffic demands are set as 580, 720, 860 veh/h to represent different levels in general lane.

When demand is low (580 veh/h), travel time changes minimally with increasing CAV MPRs. This is because CAVs do not need frequent access to the bus lane to pass the stop bar smoothly. However, under high demand (860 veh/h), when MPR reaches 20%, the average travel time decreases by ~14% in Scenario A and ~12% in Scenario B. The reduction peaks at ~31% and ~16%, respectively, when MPR reaches 40%. As MPR increases further, travel times stabilize and converge with those observed under low-demand conditions, indicating that the DBPL strategy is effective in handling high traffic demand.

To explore why the DBPL strategy continues to enhance system efficiency at high MPRs, Fig. 17 compares vehicle trajectories under DBPL and EBL strategies with 100% MPR and 860 veh/h demand. As seen, the DBPL strategy reduces congestion in the general lane by redirecting some CAVs to the bus lane. These CAVs pass the stop bar earlier than if they remained in the general lane. Furthermore, Fig. 17 (c) and (d) show that even under high demand, bus priority is not compromised by CAVs using the bus lane.

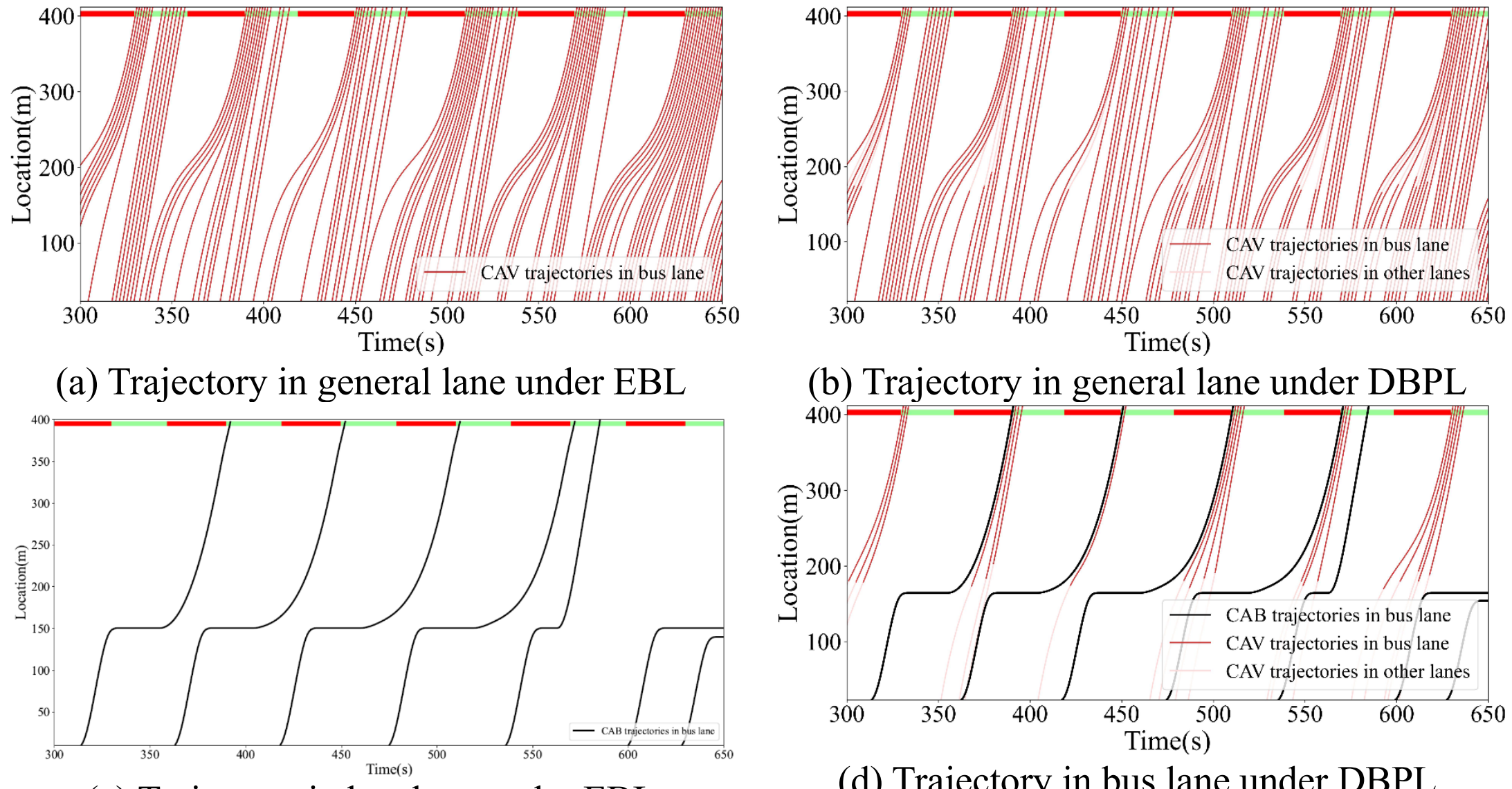

(a) Trajectory in general lane under EBL (b) Trajectory in general lane under DBPL

(c) Trajectory in bus lane under EBL (d) Trajectory in bus lane under DBPL

Fig. 17 Vehicle trajectories between DBPL and EBL in scenario A under high demand

### *4.3.2 Sensitivity to bus arrival intervals*

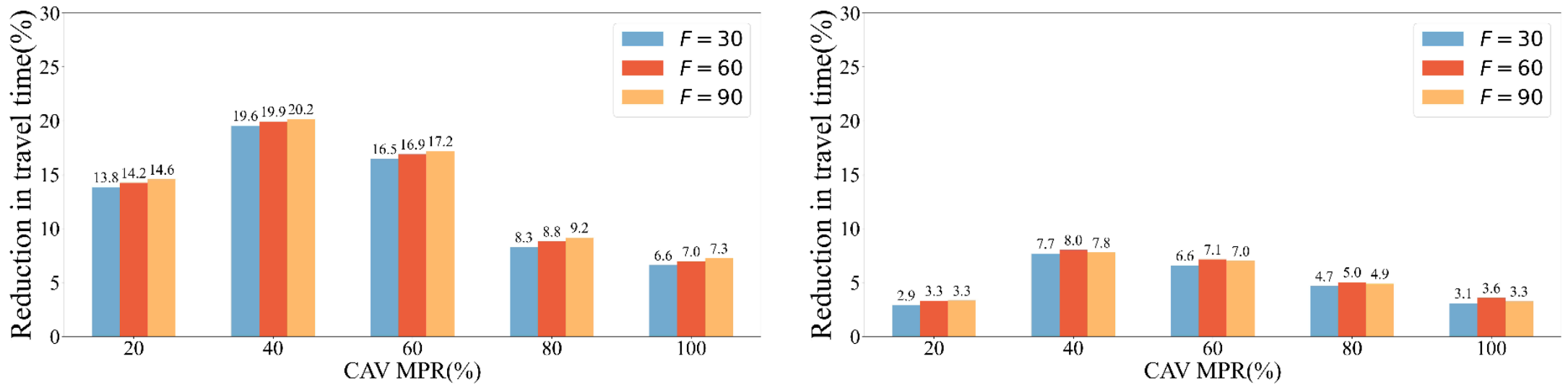

(a) Travel time reduction in scenario A (b) Travel time reduction in scenario B

Fig. 18 Impacts on travel time reduction of the CAV MPRs at different bus arrival interval

Fig. 18 evaluates how changes in bus arrival intervals (F) impact travel time reduction. Bus intervals are set at 30, 60, and 90 seconds, with a variance of 20, to represent varying levels of bus frequency. $Q_{veh}$is fixed at 720 veh/h, and $x_s$ is 150 m. Results show that bus arrival intervals have minimal impact on the travel time reduction achieved by the DBPL strategy. Even when buses arrive frequently (F=30s), private car travel times are reduced by up to ~19.6% in Scenario A and ~7.7% in Scenario B.

To further illustrate the benefits of the DBPL strategy, Fig. 19 presents vehicle trajectories with 860 veh/h demand, 60% CAV MPR, and different bus intervals. As seen in Fig. 19 (a) and (b), even with a high bus frequency (F = 30s), bus trajectories under the DBPL strategy remain nearly identical to those under the EBL strategy. As bus arrival intervals increase (F=60s to 90s), more spatial-temporal resources become available, but the number of CAVs utilizing the bus lane does not increase significantly. This is because allowing additional CAVs into the bus lane would not necessarily lead to higher traffic efficiency and could even reduce overall efficiency by diminishing the leading effect of CAVs in the general lane.

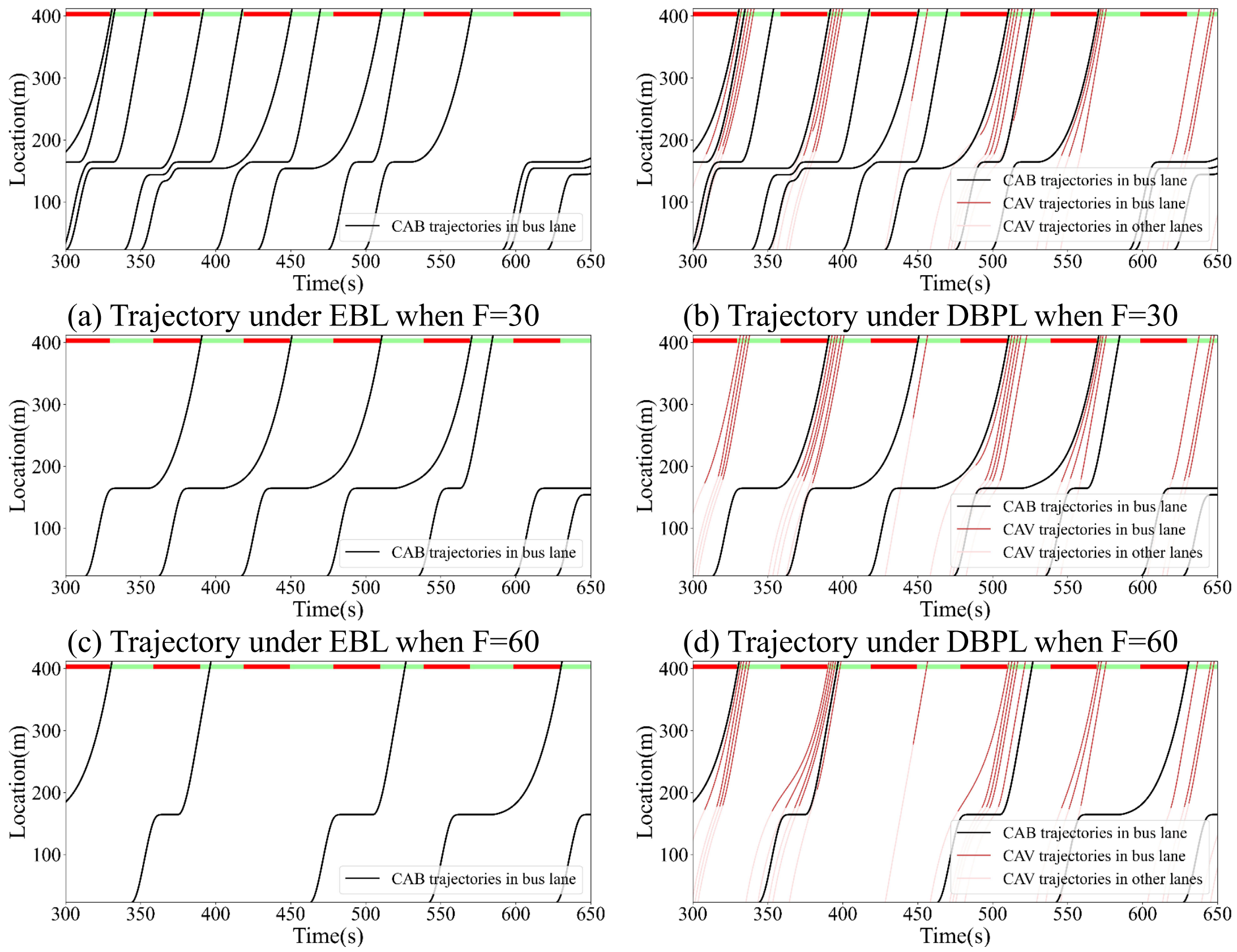

(a) Trajectory under EBL when F=30 (b) Trajectory under DBPL when F=30

(c) Trajectory under EBL when F=60 (d) Trajectory under DBPL when F=60

(e) Trajectory under EBL when F=90 (f) Trajectory under DBPL when F=90

Fig. 19 Vehicle trajectories in bus lane between DBPL and EBL in scenario A under different bus arrival intervals

### *4.3.3 Sensitivity to bus stop position*

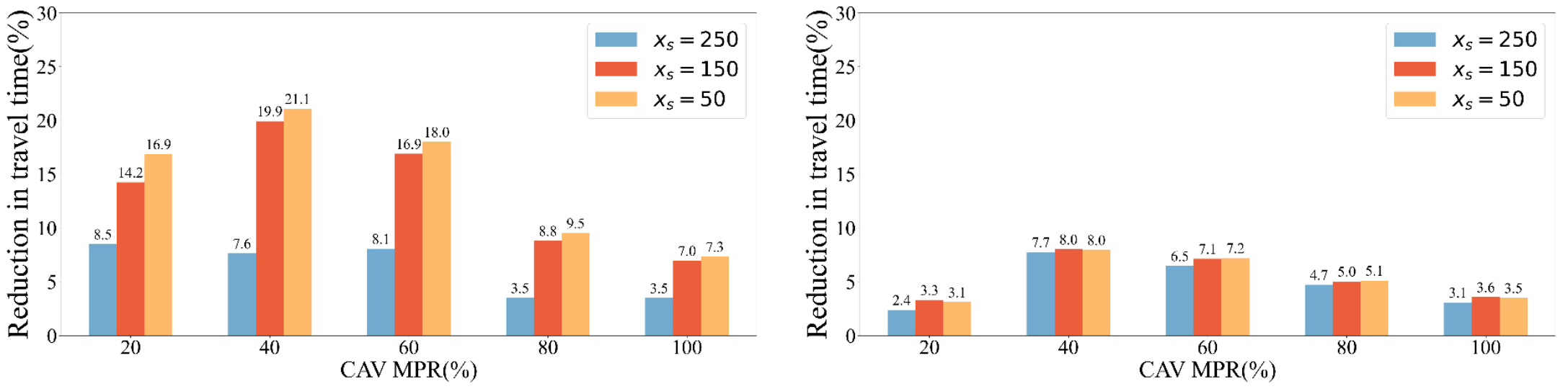

(a) Travel time reduction in scenario A (b) Travel time reduction in scenario B

Fig. 20 Impacts on travel time reduction of the CAV MPRs at different bus stop position

Fig. 20 analyzes how varying bus stop positions affect travel time reduction. Simulations are conducted with traffic demand at 720 veh/h, MPRs ranging from 20% to 100%, and bus arrival intervals fixed at 60 seconds. Three $x_s$ are tested: 50m, 150m, and 250m.

In Scenario A, as the bus stop moves farther from the intersection, the reduction in private car travel time becomes more significant. This is illustrated in Fig. 21 (a) and (b): the increased distance creates more available spatial-temporal resources in the bus lane, allowing a greater number of CAVs to use it, which helps alleviate congestion in the general lane. Additionally, comparing Fig.

21 (a) with (c) or (b) with (d) shows that the DBPL strategy has minimal impact on bus operations.

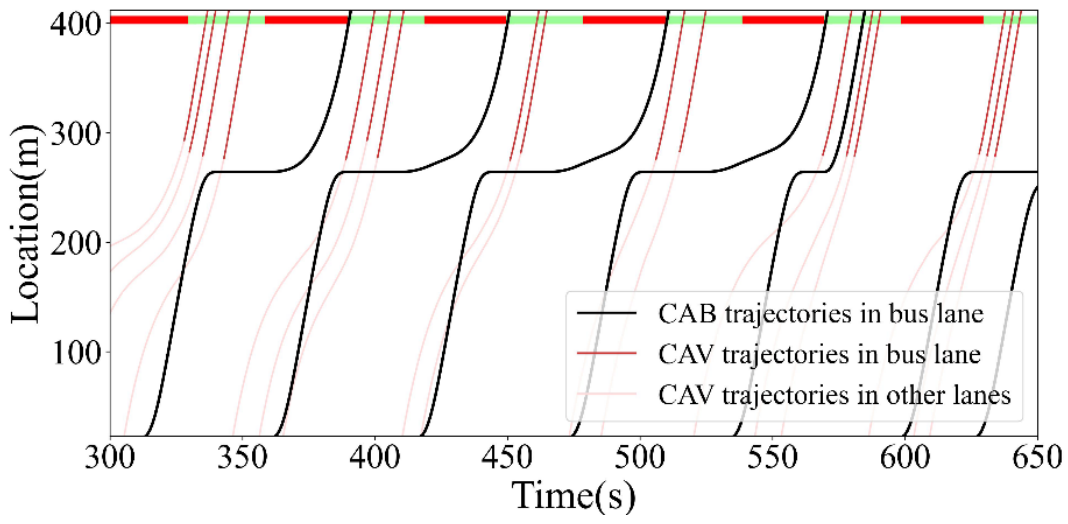

(a) Trajectory under DBPL when $x_s = 250m$

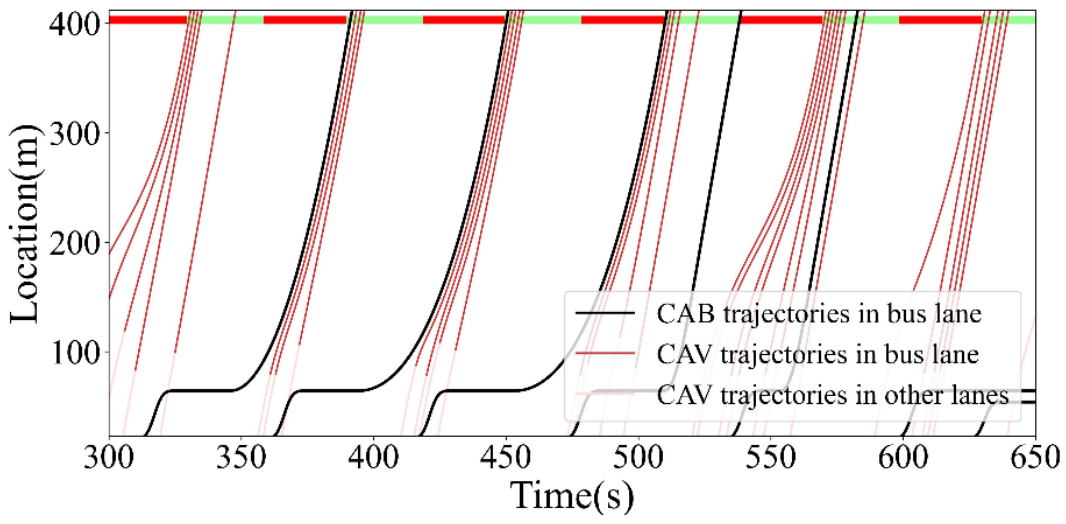

(b) Trajectory under DBPL when $x_s = 50m$

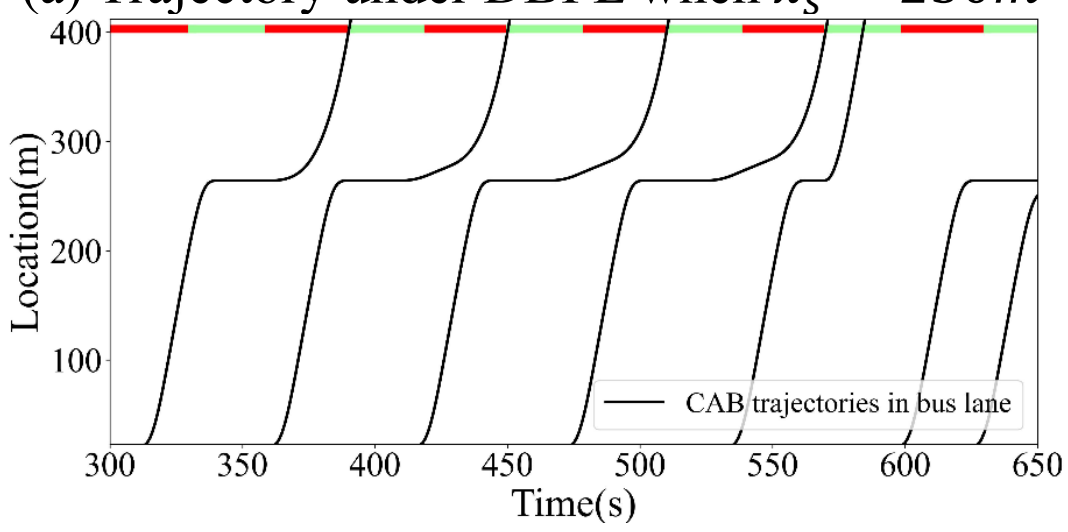

(c) Trajectory under EBL when $x_s = 250m$

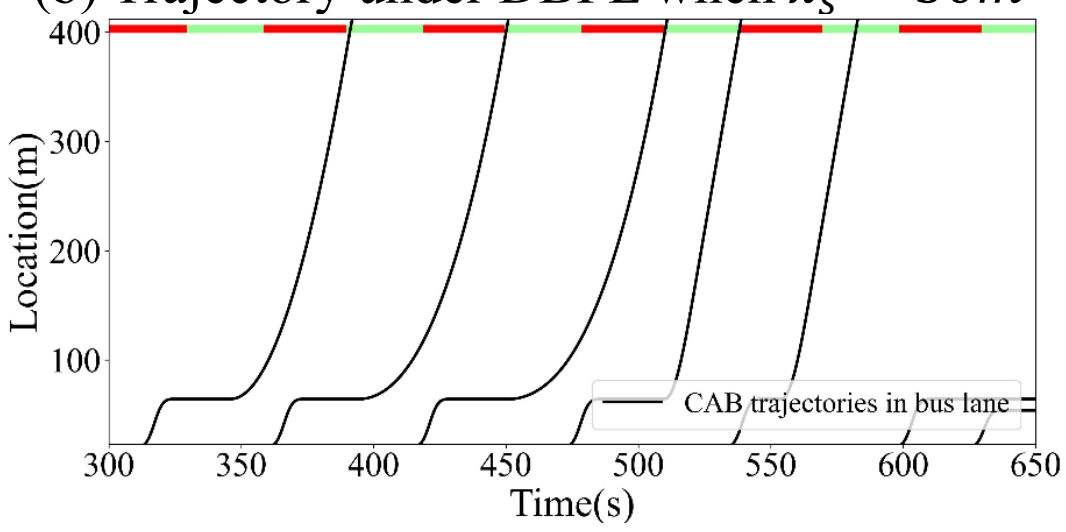

(d) Trajectory under EBL when $x_s = 50m$

Fig. 21 Vehicle trajectories in bus lane between DBPL and EBL in scenario A under different bus stop position when CAV MPR is 60%

In contrast, under the same conditions, the DBPL strategy in Scenario B is less affected by the bus stop location. This is because traffic in Scenario B is already smooth even before implementing the DBPL strategy, meaning fewer CAVs need to enter the bus lane to maintain an efficient flow. This behavior is further illustrated in Fig. 22 (a) and (b): even with additional available space in the bus lane, the number of CAVs using it remains unchanged. This shows that the DBPL strategy adaptively considers overall traffic efficiency based on current traffic conditions, avoiding unnecessary CAV entries into the bus lane merely due to increased spatial-temporal availability, which demonstrates its reliability.

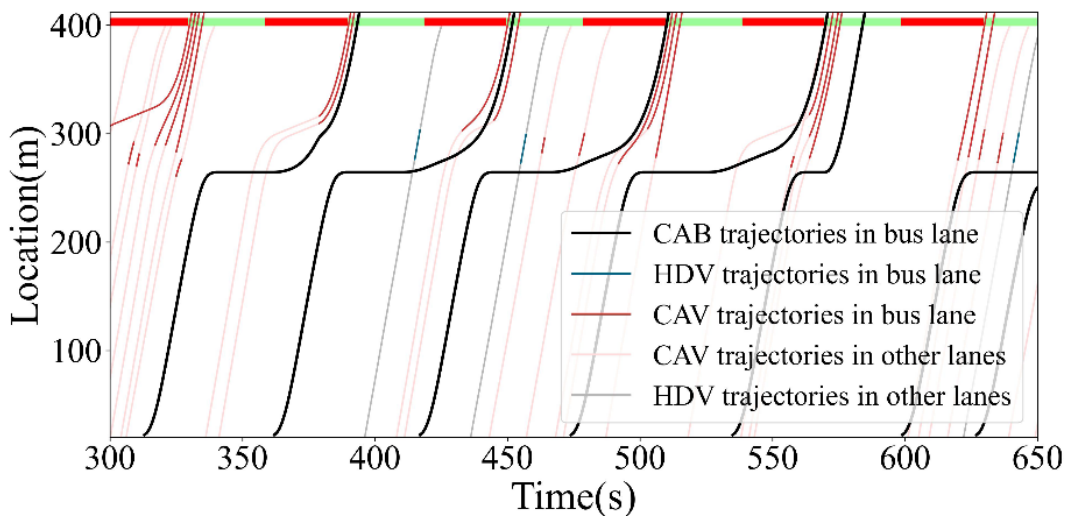

(a) Trajectory under DBPL when $x_s = 250m$

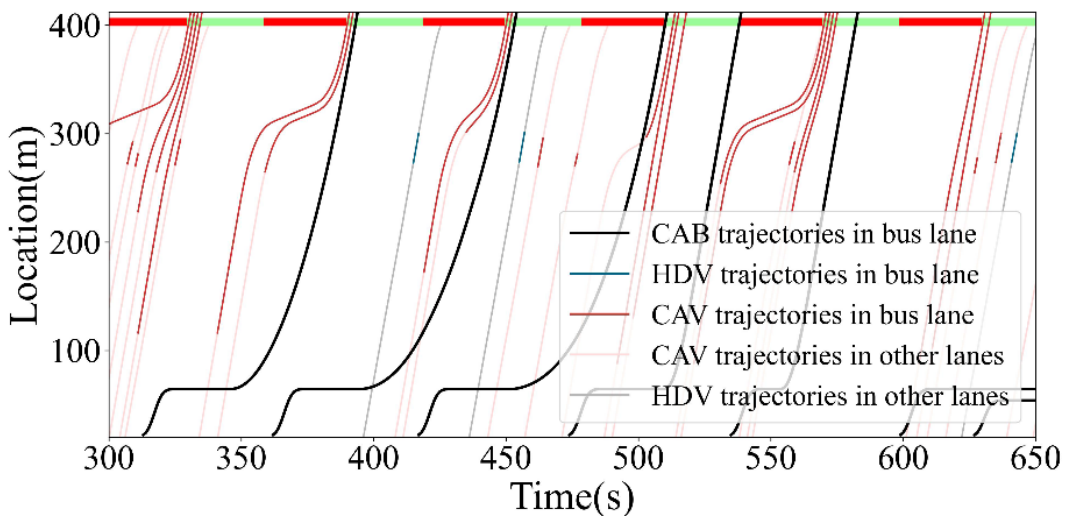

(b) Trajectory under when $x_s = 50m$

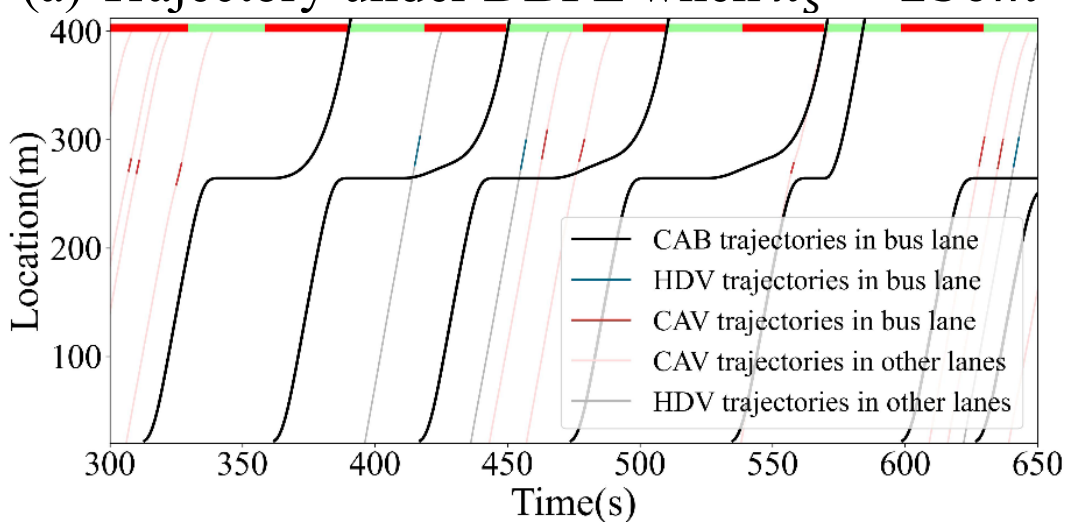

(c) Trajectory under EBL when $x_s = 250m$

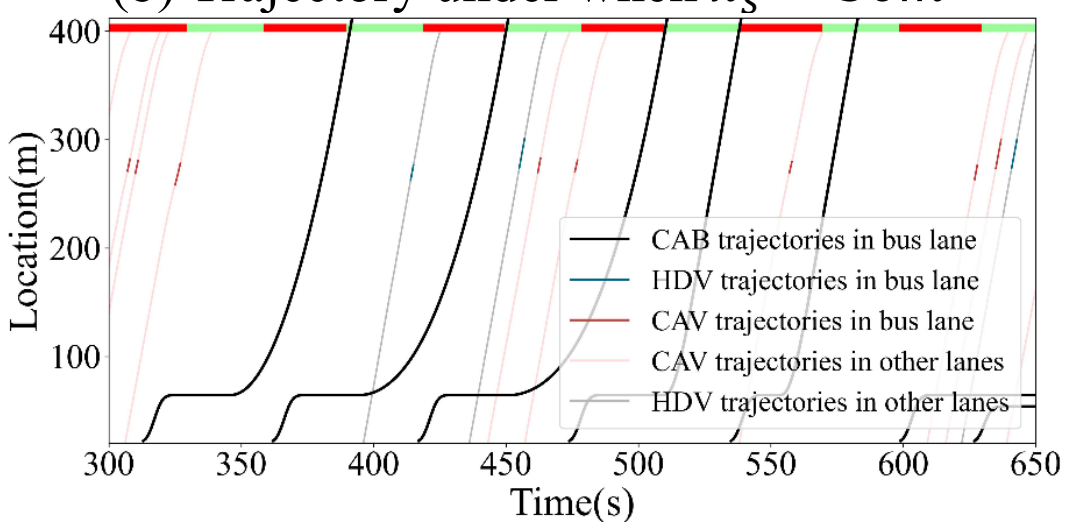

(d) Trajectory under EBL when $x_s = 50m$

Fig. 22 Vehicle trajectories in bus lane between DBPL and EBL in scenario B under different bus

stop position when CAV MPR is 60%

### *4.3.4 Sensitivity to right-turn ratio*

In this section, we examine how different right-turn ratios affect the DBPL strategy in scenario B. The simulation is conducted with a traffic demand of 720 veh/h, MPR of 40%, and bus arrival intervals of 60 seconds. Fig. 23 (a) and (b) present the travel times and reductions for right-turning and non-right-turning vehicles.

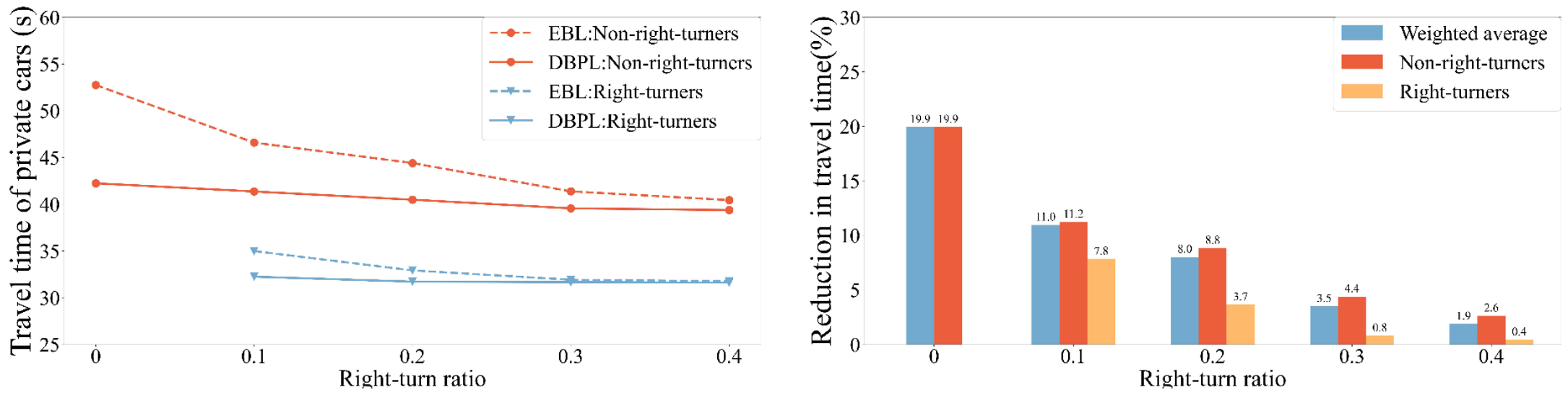

(a) Travel time in scenario B (b) Travel time reduction in scenario B

Fig. 23 Travel time and reduction of private cars under different right-turn ratio

As shown in Fig. 23 (a), as the right-turn ratio increases, the travel times of both right-turning and non-right-turning vehicles decrease significantly under both EBL and DBPL strategies. This reduction occurs because the higher proportion of right-turning vehicles, which are not controlled by traffic signals, reduces the number of through vehicles, thereby lowering the overall average travel time. The reduction trend in Fig. 23 (b) aligns with the observations in Fig. 23 (a). Even with a high right-turn ratio of 0.4, where right-turn traffic demand is substantial, the DBPL strategy maintains smooth traffic flow. Notably, the travel time of non-right-turning vehicles still decreases by approximately 3%, demonstrating the continued effectiveness of the DBPL strategy in optimizing overall traffic performance.

## 5. Conclusions

This study presents a novel DBL control strategy tailored for mixed manual and automated traffic environment. The primary objective is to alleviate congestion in general lanes while utilizing the spare capacity of bus lanes to benefit CAVs. To achieve this, we developed a MILP model to optimize the ROW allocation for CAVs within a bus priority framework. The model selectively permits eligible CAVs to access bus lanes, aiming to minimize total road delays.

The strategy reduces delays through two key mechanisms: (1) Strategically retaining some CAVs in general lanes to leverage their leading effect in smoothing traffic flow. (2) Enabling some CAVs to enter bus lanes to reduce the vehicle load in general lanes and facilitate CAVs proceeding intersections more rapidly. While CAVs and CABs perform trajectory planning autonomously, the traffic controller assigns ROW without directly altering their planned trajectories. To adapt to changing traffic conditions, a dynamic control framework incorporating a rolling horizon procedure was employed to update the optimal ROW allocation in real-time.

The proposed strategy was evaluated across two road scenarios. Results show that the DBPL control strategy reduced private car travel time by up to 20% while maintaining smooth bus operations. Sensitivity analysis further demonstrated its effectiveness compared to the EBL strategy under varying traffic demands, with private car travel time decreased by up to 31% as traffic volume increased. Even under varying bus arrival intervals, the DBPL strategy effectively utilized available bus lane capacity without disrupting bus priority, reducing private car travel time by 3% to 20%. The system also ensures that only the necessary number of CAVs enter the bus lane, optimizing overall traffic flow. Additionally, DBPL performed well across different bus stop locations and right-turn ratios, underscoring its adaptability to diverse traffic conditions.

It is worth to point out that this study does not consider the potential benefits of optimizing traffic signal timing alongside ROW allocation to further enhance road efficiency. Future research could explore the co-optimization of signal timing and vehicle ROW allocation. Additionally, it would be valuable to investigate ROW optimization strategies at the corridor or network level under partially connected vehicle environments.

## Appendix A

Main notations used in this paper are summarized as follows:

Table 3 Summary of key notations

| **General notations** | |
|---|---|
| $i$ | Vehicle index |
| $\bar{i}$ | The preceding vehicle of vehicle $i$ |
| $\bar{i}'$ | The nearest through vehicle ahead of vehicle $i$ |
| $r_i$ | The virtual vehicle that corresponds to the real vehicle $i$ |
| $\overline{r_i}$ | The preceding vehicle of vehicle $i$ in the target lane |
| $\underline{r_i}$ | The following vehicle of vehicle $i$ in the target lane |
| $p_i$ | The anticipated preceding vehicle in the bus lane for vehicle $i$ granted ROW |
| **Parameters** | |
| $k_0$ | Start time step |
| $k_c$ | Lane change time communicated by the control center to vehicles |
| $\Delta k$ | Length of each time step |
| $h$ | Length of the planning horizon |
| $\omega_p$ | Weighting parameter balancing the priority between buses and cars |
| $x_c$ | Position of the stop bar |
| $x_n$ | Start position of the no-changing zone |
| $x_w$ | Position of the right-turn pocket entrance |
| $x_s$ | Position of the bus stop |
| $v_i^{max}$ | Maximum speed of vehicle $i$ |
| $a_i^{max}$ | Maximum acceleration of vehicle $i$ |
| $a_i^{min}$ | Maximum deceleration of vehicle $i$ |
| $d_{safe}$ | Safe distance needed between vehicles for lane-changing |
| $k_{lc}$ | Lane-changing time step |
| $t_c$ | Signal cycle time (including green $t_g$, red and amber phases $t_r$) |
| $t_{bl}$ | Passing time of the last real vehicle in the bus lane before $k_0$ |
| $\underline{t}_{bl}$ | Time the last real vehicle passed $x_w$ before $k_0$ |
| $\underline{v}_{bl}$ | Speed the last real vehicle passed $x_w$ before $k_0$ |
| $d_i$ | Spatial displacement of vehicle $i$ in Newell's car-following model |
| $\tau_i$ | Temporal displacement of vehicle $i$ in Newell's car-following model |
| $\hat{\tau}_i$ | Minimum car-following time headway for vehicle $i$ |
| $\tau_r$ | Reaction time to the green light |
| $\tau_a$ | Time required to start and accelerate to pass the stop bar |
| $l_i$ | Length of vehicle $i$ |
| $l_v$ | Length of private cars |
| $l_B$ | Length of buses |
| $v_i(k)$ | Speed of vehicle $i$ at time $k$ |
| $x_i(k)$ | Position of vehicle $i$ at time $k$ |
| $\vartheta_i(k)$ | Indicator for lane-changing opportunities for vehicle $i$ |
| $t_i^{tra}(k)$ | Time vehicle $i$ passes the stop bar (from shared trajectory data at time $k$) |
| $t_i^{p1}(k)$ | Earliest time vehicle $i$ can pass the stop bar without considering preceding vehicles |
| $\underline{t}_i^{tra}(k)$ | Time vehicle $i$ passed $x_w$ (from shared trajectory data at time $k$) |
| $\underline{v}_i^{tra}(k)$ | Speed of vehicle $i$ passing $x_w$(from shared trajectory data at time $k$) |
| **Sets** | |
| $K$ | Set of all time steps within the planning horizon: $K = \{k_0, k_0 + \Delta k, k_0 + 2\Delta k, \dots, k_0 + h\}$ |
| $\mathcal{I}_a$ | Set of CAVs |
| $\mathcal{I}_h$ | Set of HDVs |
| $\mathcal{I}_b$ | Set of buses |
| $\mathcal{I}_{gl}$ | Set of vehicles in general lane |
| $\mathcal{I}_{bl}$ | Set of vehicles in the bus lane, with subsets of real vehicles and virtual vehicles denoted as $\mathcal{I}_{bl}^{real}$ and $\mathcal{I}_{bl}^{vir}$,respectively: $\mathcal{I}_{bl} = \mathcal{I}_{bl}^{vir} \cup \mathcal{I}_{bl}^{real}$ |
| $\mathcal{I}_c$ | Set of vehicles granted ROW |
| $I'_{gl}(k)$ | Vehicles in the general lane that have passed the stop bar at time step $k$ |
| $\bar{I}_{gl}(k)$ | Vehicles in the general lane yet to pass the stop bar at time step $k$ |
| $I'_{bl}(k)$ | Vehicles in the bus lane that have passed the stop bar at time step $k$ |
| $\bar{I}_{bl}(k)$ | Vehicles in the bus lane that have passed $x_w$ at time step $k$ |
| $\underline{I}_{bl}(k)$ | Vehicles in the bus lane yet to pass $x_w$ at time step $k$ |
| $\underline{I}_{bl}^T(k)$ | Subset of $\underline{I}_{bl}(k)$ consisting of through vehicle |
| **Variables** | |
| $T_c$ | Average travel time of cars |
| $T_b$ | Average travel time of buses |
| $t_i^{dep}(k)$ | Time vehicle $i$ passes the stop bar |
| $t_i^{p2}(k)$ | Earliest time vehicle $i$ can pass the stop bar without signal control |
| $t_i^{p3}(k)$ | Departure time of vehicle $i$ considering signal effects |
| $t_i^{pre}(k)$ | Time when the preceding vehicle of $i$ passes the stop bar |
| $s_i^{acc}(k)$ | Distance required for vehicle $i$ to accelerate to maximum speed |
| $t_i^{acc}(k)$ | Time required for vehicle $i$ to reach maximum speed |
| $t_i^R(k)$ | Start time of the red phase within the signal cycle for vehicle $i$ |
| $\underline{t}_i^{dep}(k)$ | Time vehicle $i$ passes $x_w$ |
| $\underline{t}_i^{p1}(k)$ | Earliest time vehicle $i$ passes $x_w$ without considering preceding vehicles |
| $\underline{t}_i^{p2}(k)$ | Time vehicle $i$ passes $x_w$ considering preceding vehicles |
| $\underline{t}_i^{pre}(k)$ | Time when the preceding vehicle of $i$ passes $x_w$ |
| $\underline{v}_i^{dep}(k)$ | Speed of vehicle $i$ passing $x_w$ |
| $\underline{v}_i^{pre}(k)$ | Speed of the preceding vehicle of $i$ passing $x_w$ |
| $\underline{v}_i^{p1}(k)$ | Speed at which vehicle $i$ passes $x_w$ without considering preceding vehicles |
| $\underline{v}_i^{p2}(k)$ | Speed at which vehicle $i$ passes $x_w$ considering preceding vehicles |
| $\underline{v}_i^{p3}(k)$ | Speed of vehicle $i$ passing $x_w$ without considering virtual vehicles |
| $\phi_i(k)$ | Binary variable indicating if vehicle $i$ changes to the bus lane at step $k$ |
| $\lambda(k)$ | Binary variable indicating if a lane change command is issued at time step $k$ |